\documentclass[12pt]{article}
\usepackage{amsfonts}
\usepackage{amssymb}
\usepackage{epsfig}
\catcode`\@=11
\newcount\hour
\newcount\minute
\newtoks\amorpm \hour=\time\divide\hour by 60\minute
=\time{\multiply\hour by 60 \global\advance\minute by-\hour}
\edef\standardtime{{\ifnum\hour<12 \global\amorpm={am}%
        \else\global\amorpm={pm}\advance\hour by-12 \fi
        \ifnum\hour=0 \hour=12 \fi
        \number\hour:\ifnum\minute<10
        0\fi\number\minute\the\amorpm}}
\edef\militarytime{\number\hour:\ifnum\minute<10
0\fi\number\minute}
\def\draftlabel#1{{\@bsphack\if@filesw {\let\thepage\relax
   \xdef\@gtempa{\write\@auxout{\string
      \newlabel{#1}{{\@currentlabel}{\thepage}}}}}\@gtempa
   \if@nobreak \ifvmode\nobreak\fi\fi\fi\@esphack}
        \gdef\@eqnlabel{#1}}
\def\@eqnlabel{}
\def\@vacuum{}
\def\marginnote#1{}
\def\draftmarginnote#1{\marginpar{\raggedright\scriptsize\tt#1}}
\def\draft{
        \pagestyle{plain}
        \overfullrule=2pt
        \oddsidemargin -.1truein
        \def\@oddhead{\sl \phantom{\today\quad\militarytime} \hfil
        \smash{\Large\sl DRAFT} \hfil \today\quad\militarytime}
        \let\@evenhead\@oddhead
        \let\label=\draftlabel
        \let\marginnote=\draftmarginnote
        \def\ps@empty{\let\@mkboth\@gobbletwo
        \def\@oddfoot{\hfil \smash{\Large\sl DRAFT} \hfil}
        \let\@evenfoot\@oddhead}
        \def\@eqnnum{(\theequation)\rlap{\kern\marginparsep\tt\@eqnlabel}%
        \global\let\@eqnlabel\@vacuum}  }
\newcommand\euro{{
C%
    \makebox[0pt][l]{\kern-.70em\mbox{--}}%
    \makebox[0pt][l]{\kern-.68em\raisebox{.25ex}{--}}}}

\renewcommand{\theequation}{\thesection.\arabic{equation}}
\renewcommand{\thefootnote}{\fnsymbol{footnote}}
\newcommand{\newsection}{    
\setcounter{equation}{0}\section}
\def\appendix#1{\addtocounter{section}{1}\setcounter{equation}{0}
\renewcommand{\thesection}{\Alph{section}}
\section*{Appendix \thesection\protect\indent \parbox[t]{11.15cm}{#1}}
\addcontentsline{toc}{section}{Appendix \thesection\ \ \ #1}}

\def \lc {{light-cone}}

\def \la {\label}

\def\be{\begin{equation}}
\def\ee{\end{equation}}

\def\eu {{\cal C}}
 \def \lc {light-cone\ }

\catcode`\@=11
\def \lc {light cone\ }

\def\bea{\begin{eqnarray}}
\def\eea{\end{eqnarray}}
\def\beann{\begin{eqnarray*}}
\def\eeann{\end{eqnarray*}}
\def\beq{\begin{equation}}
\def\eeq{\end{equation}}
\def\ba{\begin{array}}
\def\ea{\end{array}}
\def\ben{\begin{enumerate}}
\def\een{\end{enumerate}}

 \def \la {\label}
 \def\be{\begin{equation}}
\def\ee{\end{equation}}

\def \la {\label}

\font\mybb=msbm10 at 11pt

\def\bb#1{\hbox{\mybb#1}}

\def\bZ {\bb{Z}}
\def\bR {\bb{R}}

\def\bC {\bb{C}}

\def\e  {\epsilon}

\def \ee {\epsilon}

\def\lc{\lrcorner}
\def\rc{\bar \lrcorner}

\def\be{\begin{equation}}
\def\ee{\end{equation}}

\def \la{\label}

\newtheorem{corollary}{Corollary}
\newtheorem{proposition}{Proposition}
\newtheorem{theorem}{Theorem}
\newtheorem{definition}{Definition}

\newcommand{\bproof}{\noindent{\it Proof: }}
\newcommand{\eproof}{\  q.~e.~d. \vspace{0.2in}}

\begin{document}
\date{November 2002}
\begin{titlepage}
\begin{center}
\leftline{}
\vskip 1cm
\vspace{2.0cm}
{\Large \bf  Spin Cohomology}\\[.2cm]

\vspace{1.5cm}
 {\large   George Papadopoulos}

 \vspace{0.5cm}

 Department of Mathematics\\
 King's College London\\
 Strand\\
 London WC2R 2LS
\end{center}

\vskip 1.5 cm
\begin{abstract}
We explore differential and algebraic operations
on the exterior product  of spinor representations and their twists
that give rise to cohomology, the spin cohomology.
A linear differential operator $d$ is introduced which is associated to a connection $\nabla$
and a parallel spinor $\zeta$, $\nabla\zeta=0$,
and the algebraic operators $D_{(p)}$ are constructed from skew-products
 of $p$ gamma matrices.  We exhibit a large number of spin cohomology
operators and we investigate the spin cohomologies associated with connections whose holonomy
is a subgroup of $SU(m)$, $G_2$, $Spin(7)$ and $Sp(2)$. In the $SU(m)$ case,  we find
that the spin cohomology of complex spin and spin$_c$ manifolds is related to
a twisted Dolbeault cohomology.  On Calabi-Yau type of  manifolds of
dimension $8k+6$, a spin  cohomology can be defined on a twisted
complex with   operator $d+D$ which is the sum of a differential and algebraic one.
We compute this cohomology  on six-dimensional
Calabi-Yau manifolds using a spectral sequence.
In the  $G_2$ and $Spin(7)$ cases, the spin cohomology is related
to the de Rham cohomology.

\end{abstract}
\end{titlepage}
\newpage
\setcounter{page}{1}
\renewcommand{\thefootnote}{\arabic{footnote}}
\setcounter{footnote}{0}

\setcounter{section}{0}
\setcounter{subsection}{0}
\newsection{Introduction}

 On spin manifolds apart from the exterior derivative $d$
  and the associated de Rham complex $(\Lambda^*(M),
d)$, one can define the Dirac operator $(\Delta(M), D)$,
where $\Delta(M)$ is the spin bundle\footnote{We adopt
the notation to denote a representation and its
 associated bundle with the same symbol, e.g.
$\Delta=\Delta(\bR^n)$ denotes the spin representation
of $Spin(n)$ and $\Delta=\Delta(M)$ denotes also the spin
bundle over $M$. In addition we shall denote the bundles
 and their sections with the same symbol. $\Lambda^*$
denotes the space of forms}. On complex manifolds the Dirac
operator decomposes as $D={\cal D}+\bar {\cal D}$ and
the spin representation can be graded such that
 $(\Delta(M), \bar {\cal D})$ can turn into a (graded) complex.
The associated cohomology is called {\it spinor} cohomology \cite{lawson}.

In even dimensions the complex Dirac spin represenation is
reducible and  decomposes as $\Delta=\Delta^-\oplus \Delta^+$.
 Apart from the spin representation $\Delta$ both
the exterior power,
$\eu=\Lambda^*(\Delta^*)$ ($\eu_\pm=\Lambda^*(\Delta_\pm)$), and the symmetric product, ${\rm Sym}^*(\Delta^*)$
(${\rm Sym}^*(\Delta_\pm)$),  of the dual spinor representation $\Delta^*$ ($\Delta_\pm= (\Delta^\pm)^*$)
have found applications in various problems in physics.
The former has applications in supermanifolds. In particular all real supermanifolds
that appear in the context of supersymmetry
are isomorphic to $\eu$ ($\eu_\pm$) \cite{bat}.  The latter is a model for the odd forms on supermanifolds and
  has
appeared in the context of string theory \cite{berk} and  the theory of deformations of
the field equations of supersymmetric gauge theories and supergravity in superspace \cite{tsimpis, howe}.
It turns out that the theory of deforming the field equations of the supersymmetric
gauge theories and supergravity can turn into a cohomological problem on
$\Lambda^*\otimes {\rm Sym}^*$
for the so called  {\it spinorial} cohomology.

Motivated by these developments in physics, the aim of this paper is to
investigate various cohomology operators that can be
defined on $\eu$, $\eu_\pm$ and its various twistings.
Let $(M,g)$ be a spin manifold equipped with a spin connection
$\nabla$, which is {\it not}
necessarily the Levi-Civita connection of the metric $g$.
One can define a linear differential (spin) operator on $\eu(M)$ or $\eu_\pm(M)$ as
\be
d\phi=\zeta\Gamma^\mu\bar\wedge\nabla_\mu\phi~,~~~~~\phi\in \eu(M)~,
\la{ldif}
\ee
where $\bar\wedge$ is the wedge product in $\eu$, $\zeta$ is a cospinor
and $\{\Gamma^\mu: \mu=1,\dots, {\rm dim} M\}$ are the gamma matrices.
In many applications $\zeta$ is taken
to be a parallel cospinor with respect to $\nabla$, $\nabla\zeta=0$.
As we shall see there are various cohomology theories
that can be defined depending on the choice of spinor $\zeta$, the connection
$\nabla$. The operator $d$ can always be defined on $\eu$. However, the restriction
of $d$ on $\eu_\pm$ depends on the choice of $\zeta$ and the properties of
the spinor inner products which in turn depend on the dimension of the manifold $M$.
One of our  aims is  to investigate the conditions for $d$ to be nilpotent, $d^2=0$.
These conditions can be expressed in terms of restrictions on the cospinor $\zeta$
and on the curvature $R$ of the connection $\nabla$, in addition to $\nabla\zeta=0$ .

In addition to differential operators, we shall present a large number
of algebraic cohomology operators $D_{(p)}$ on
some twisted complexes like for example $\Lambda^*(M)\otimes\eu(M)$ and $\Lambda^*(M)\otimes\eu_\pm(M)$.
Some of these are constructed from skew-products of $p$ gamma matrices. We investigate
the conditions  for $D_{(p)}^2=0$ and relate them into the symmetry properties
of gamma matrices. The latter again depend on the dimension of the manifold $M$.
In addition we shall show that $D_{(p)} d+d D_{(p)}=0$ and
so the cohomology of $(\Lambda^*\otimes \eu, d+D_{(p)})$ and  $(\Lambda^*\otimes \eu_\pm, d+D_{(p)})$
can be computed using a spectral sequence.
We shall refer collectively to all of these cohomology theories with operators $d$, $D_{(p)}$ and $d+D_{(p)}$ as
{\it spin} cohomologies.

We shall develop the general theory of spin cohomology. In particular,
we shall compute the conditions on the curvature of the underlying manifold
for $d^2=0$. We shall also explain the relation to parallel spinors.

Next, we shall focus
on a certain class of parallel spinors. In particular,
we shall consider manifolds which admit a spin  connection $\nabla$
induced from the tangent bundle
with holonomy contained in the groups  $SU(m)$ (n=2m),
$Sp(k)$ (n=4k), $Spin(7)$ (n=8) and $G_2$
(n=7), where in parenthesis is the dimension of the manifold. A special
class of examples of manifolds with spin cohomology
are those of which $\nabla$ is the Levi-Civita connection.

We shall show that for spin complex manifolds which admit a holomorphic
connection $\nabla$ with ${\rm hol}\nabla\subseteq SU(m)$, there are
two differential spin cohomologies with operators $d_1$ and $d_2$ related to two parallel
spinors of the connection $\nabla$. We refer to these spin cohomologies
as complex spin cohomologies. We shall show that $d_1$ and $d_2$ restrict on $\eu_\pm$. In particular
one can construct complexes $(\eu_+, d_1)$ and $(\eu_+, d_2)$ for ${\rm dim M}=8k+2, 8k+6$ and
complexes $(\eu_-, d_1)$ and $(\eu_-, d_2)$ for ${\rm dim M}=8k, 8k+4$.
We  give the Laplace operators
 associated with $d_1$ and $d_2$ using a $Spin(n)$-invariant inner product.
We show that the complex spin cohomology of $(\eu_-, d_2)$ in all dimensions is
related to  twisted Dolbeault cohomology. We extend this relation
between this spin cohomology and Dolbeault cohomology to complex spin$_c$ manifolds as well.
The complex spin cohomologies can be  twisted with any holomorphic vector bundle.
Apart form the differential complex spin cohomologies, there is an algebraic spin cohomology operator $D=D_{(1)}$
and the complex $\Lambda^{*,0}\otimes \eu_-$
on all such manifolds of dimension $n=8k+6$ and $d_2 D+ D d_2=0$.
The cohomology of  $(\Lambda^{*,0}\otimes \eu_-, d_2+D)$ can be computed using a spectral
sequence. As an example, we computed the cohomology of $(\Lambda^{*,0}\otimes \eu_-, d_2+D)$
 on six-dimensional Calabi-Yau manifolds.

On  manifolds which admit a connection with holonomy $Sp(k)$, there are $k+1$ differential spin
operators associated to $k+1$ parallel spinors.
Two of these are the same as those
of the $SU(2k)$ manifolds investigated above.
 We shall not present a full analysis in this case but
we shall express a third spin differential operator
on hyperK\"ahler manifolds in terms of a Dolbeault operator.

On  manifolds which admit a connection with holonomy  $Spin(7)$, there is one differential spin operator $d$
associated to one parallel spinor and a real complex $(\eu_{\bR}, d)$. In addition,
 $d^2=0$ provided the connection $\nabla$ is the Levi-Civita
connection of a $Spin(7)$ metric. The spin cohomology is isomorphic
to de Rham cohomology.

On  manifolds which admit a connection with holonomy  $G_2$, there is again one differential spin operator $d$
and a real complex $(\eu_{\bR}, d)$.
In addition,  $d^2=0$ provided the connection
 $\nabla$ is the Levi-Civita
connection of a $G_2$ metric. The spin cohomology of $(\eu_{\bR}, d)$ is isomorhic to two copies
of the de Rham cohomology relatively shifted by one degree.

This paper has been organized as follows: In section two, we summarize the properties
of Clifford algebras and spin representations which we use later.
In section three, we explore the general properties of the
linear differential operators (\ref{ldif}),
define the twisted complexes and present the algebraic cohomology operators.
In section four, we investigate the properties of complex  spin cohomology and derive the
conditions for $d^2=0$. In addition we compute the Laplace operators. In section five, we investigate
various kinds of twisted complex cohomology.  In section six, we relate the
complex spin cohomology to the Dolbeault cohomology for spin  and spin$_c$ manifolds.
In section seven, we compute the complex spin cohomology and a twisted spin cohomology
on a six-dimensional Calabi-Yau manifold. In section eight, we investigate the spin cohomology
of manifolds that admit a connection with holonomy contained in $Sp(k)$. In section nine,
we explore the properties of some real spin cohomologies. In sections ten and eleven, we
investigate the spin cohomology of manifolds that admit a connection with holonomy
$Spin(7)$ and $G_2$, respectively.

\newsection{Preliminaries}

The investigation of spin cohomology involves a detailed description
of spinor representations. Because of this and to establish notation,
we shall review some aspects of spinor representations in various dimensions \cite{wang, harvey}.
We shall focus on the manifolds with Euclidean signature but the
analysis can be easily extended to other signatures.

Let $V=\bR^n$ be a real vector space equipped
with the standard Euclidean inner product.
 If $n=2m$ even, the
complex spin (Dirac) representation of $Spin(2m)$, $\Delta=\Delta(V)$, is
reducible and decomposes to two irreducible representations,
$\Delta=\Delta^+\oplus \Delta^-$.
To construct these spin representations  let
$e_1, \dots, e_n$ be an orthonormal basis in $\bR^n$, $n=2m$, and $J$ be
a complex structure in $V$, $J(e_i)=e_{i+m}$. We identify $V$ and its dual
using the Euclidean inner product.  Next consider
the subspace $U=\bR^m$ generated by $e_1, \dots, e_m$. Clearly $V=U\oplus J(U)$.
The Euclidean inner product on $V$ can be extended to a hermitian inner
product in $V_{\bC}=V\otimes \bC$ denoted by $<,>$, ie
\be
<z^\mu e_\mu, w^\nu e_\nu>=\sum_\mu \bar z{}^\mu w^\mu~,
\la{inin}
\ee
where $\bar{z}$~ is the standard complex conjugate of $z$ in $V_{\bC}$.
The space of spinors $\Delta(V)=\Lambda^*(U_{\bC})$, where $U_{\bC}=U\otimes \bC$.
In addition, $\Delta^+=\Lambda^{\rm even}(U_{\bC})$ and
$\Delta^-=\Lambda^{\rm odd}U_{\bC}$. The spinors in $\Delta^+$ are called chiral
while those in $\Delta^-$ anti-chiral.
The inner product (\ref{inin}) can be easily extended to $\Delta$ and it is called
the Dirac inner product on the space of spinors.
The generators of the Clifford algebra $e_\mu$ are represented on $\Delta$ as
\bea
\Gamma(e_i)\eta&=&e_i\cdot \eta= e_i\wedge\eta+e_i\lc \eta~,~~~~ i\leq m
\cr
\Gamma(e_{m+i})\eta&=&e_{i+m}\cdot \eta=-ie_i\wedge\eta+ie_i\lc \eta~,~~~~~i \leq m~,
\eea
where $e_i\lc$ is the adjoint of $e_i\wedge$ with respect to $<,>$.
It is convenient to denote the generators $\Gamma(e_{\mu})=\Gamma_\mu$ and they are often
called gamma
matrices. Clearly $\Gamma_\mu: \Delta^\pm\rightarrow \Delta^\mp$.
The linear maps $\Gamma_\mu$ are hermitian with respect
to the inner product $<,>$, $<\Gamma_\mu \eta, \theta>=<\eta, \Gamma_\mu\theta>$,
and satisfy the Clifford algebra relations
$e_\mu e_\nu+e_\nu e_\mu=\Gamma_\mu\Gamma_\nu+\Gamma_\nu \Gamma_\mu=0$, for $\mu\not=\nu$, $(e_\mu)^2=(\Gamma_\mu)^2=1$.

Next define the maps $A=\Gamma_1 \Gamma_2\dots \Gamma_m$ and
$B=\Gamma_{m+1}\dots\Gamma_n$ and the inner products on $\Delta$ as
\bea
A(\eta, \theta)&=&<A(\bar\eta), \theta>
\cr
B(\eta, \theta)&=&<B(\bar\eta), \theta>~,
\eea
which we denote with the same symbol,
where $\bar\eta$ is the standard complex conjugate of $\eta$ in $\Lambda^*(V_{\bC})$
The inner products $A,B$ are sometimes also called charge conjugation matrices.
These have the following properties:
\bea
A(\eta, \theta)&=&(-1)^{{1\over2} m (m-1)} A(\theta, \eta)
\cr
B(\eta, \theta)&=& (-1)^{{1\over2} m (m+1)} B(\theta, \eta)~.
\eea
Therefore $A$ ($B$) is symmetric for $m=4k, 4k+1$ ($m=4k, 4k+3$) and skew-symmetric
for $m=4k+2, 4k+3$ ($m=4k+2, 4k+1$).
In addition, we have
\bea
A(\Gamma_\mu\eta, \theta)&=&(-1)^{m-1} A(\eta, \Gamma_\mu\theta)~,~~1\leq \mu\leq n
\cr
B(\Gamma_\mu\eta, \theta)&=&(-1)^{m} B(\eta, \Gamma_\mu\theta)~,~~1\leq \mu\leq n
\eea
and
\bea
A(\Gamma_\mu\eta, \Gamma_\mu\theta)&=&(-1)^{m-1} A(\eta, \theta)~,~~1\leq \mu\leq n
\cr
B(\Gamma_\mu\eta, \Gamma_\mu\theta)&=&(-1)^{m} B(\eta, \theta)~,~~1\leq \mu\leq n~.
\eea
Therefore $A$ is $Pin(2m)$ invariant for $m=4k+1, 4k+3$ while $B$ is $Pin(2m)$
invariant for $m=4k, 4k+2$. Both $A,B$ are $Spin(n)$-invariant.
A consequence of the above relations is
\bea
A(\eta,  \Gamma_\mu\theta)&=&(-1)^{{1\over2} (m-1) (m+2)}
A(\theta, \Gamma_\mu\eta)
\cr
B(\eta, \Gamma_\mu\theta)&=&(-1)^{{1\over2} m (m+3)}B(\theta, \Gamma_\mu\eta)
\eea
Therefore the gamma-matrices are symmetric with respect to the
inner product $A$( $B$) for $m=4k+1, 4k+2$ ( $m=4k, 4k+1$) while
they will be skew-symmetric for $m=4k, 4k+3$ ($m=4k+2, 4k+3$).

Because of the existence of invariant non-degenerate inner products the dual
of $\Delta^*$ can be identified with $\Delta$. To make this identification
precise, let us denote with $C$ either $A$ or $B$. Given a basis
$\{\epsilon_A; A=1,\dots,
{\rm dim}\Delta\}$ in $\Delta$ let us denote with $\{\epsilon^A; A=1,\dots,
{\rm dim}\Delta\}$ the dual basis in $\Delta^*$, $\epsilon^B(\epsilon_A)=\delta^B{}_A$.
The inner product $C^{-1}$ in $\Delta^*$ induced by $C$ is
\be
\sum_{E}C^{-1}(\epsilon^E, \epsilon^A) C(\e_E, \e_B)=\delta^A_B
\ee
The
co-spinor $C(\eta)$ associated with the spinor $\eta$ under the isomorphism $C$
is defined as $C(\eta)(\theta)=C(\theta,\eta)$, $\eta,\theta\in \Delta$, ie
in the above basis $\eta_A=C_{AB} \eta^B$. The
inverse transformation
$C^{-1}$ is defined as $C^{-1}(\psi)(\chi)= C^{-1}(\psi, \chi)$,
$\psi, \chi\in \Delta^*$, ie $\psi^A= \psi_{B} (C^{-1})^{BA}$.
Notice that the maps $A,B: \Delta^\pm\rightarrow \Delta^\pm$ for $m=4k, 4k+2$
while $A,B: \Delta^\pm\rightarrow \Delta^\mp$ for $m=4k+1, 4k+3$. Therefore
in the former case the dual $\Delta_\pm$ of $\Delta^\pm$   under $A,B$ is identified
with $\Delta^\pm$, $\Delta_\pm=\Delta^\pm$ while in the latter the dual $\Delta_\pm$ of
$\Delta^\pm$ is identified with
$\Delta^\mp$, $\Delta_\pm=\Delta^\mp$.

There are two ways to construct the spin representation $\Delta$ in odd dimensions.
One is to write $V'=\bR^{2m+1}=V\oplus \bR<e_{2m+1}>$ and extend the Euclidean inner
product (\ref{inin}) from  $V$ to $V'$, $<e_{2m+1}, e_{2m+1}>=1, <V, e_{2m+1}>=0$. The gamma matrices
$\Gamma_\mu$, $ 1\leq\mu\leq 2m$, are defined as in the even-dimensional case and
\be
\Gamma_{2m+1}=i^m\Gamma_1\dots \Gamma_{2m}~.
\ee
The $Spin(2m+1)$ spin representation, $\Delta$,  is $\Delta=\Delta^+\oplus \Delta^-$, where $\Delta^+, \Delta^-$
are the $Spin(2m)$ spin representations. (There are no chiral spinors in odd dimensions.)
The invariant inner product on $Spin(2m+1)$ representation
 $\Delta$ is the $Pin(2m)$ invariant inner product on $\Delta^+\oplus \Delta^-$.

 Alternatively, we take $V=U\oplus J(U)$ as for $n=2m$ and write $U=U_0\oplus \bR<e_{2m}>$.
 Then $V_0=U\oplus J(U_0)$ has dimension $2m-1$. The gamma matrices  are
  $\tilde\Gamma_\mu= i\Gamma_\mu\Gamma_{2m}$, $1\leq\mu\leq 2m-1$, where $\Gamma_\mu$ are the gamma
 matrices of the  $Spin(2m)$ spin representation. These induce a representation of $Pin(2m-1)$
 onto the $\Delta^\pm$ representations of $Spin(2m)$.

For later convenience, we introduce the notation
$(C\Gamma_\mu )(\eta, \theta)=C(\eta, \Gamma_\mu\theta)$ and similarly
 $(C\Gamma_{\mu_1\dots\mu_p})(\eta, \theta)
=C(\eta,\Gamma_{\mu_1\dots\mu_p}\theta)$, where
\be
\Gamma_{\mu_1\dots \mu_p}= {1\over p!} \sum_{\sigma} (-1)^{|\sigma|} \Gamma_{\mu_{\sigma(1)}}
\dots \Gamma_{\mu_{\sigma(p)}}
\ee
and $\sigma$ is a permutation. The symmetry of the inner product and that of the gamma
matrices can be re-expressed as
 $C(\eta, \theta)=(-1)^{s_C} C(\theta,\eta)$, where $s_C=0$ if $C$ is symmetric and $s_C=1$
 if $C$ is skew-symmetric,  and similarly
$C\Gamma_\mu (\eta, \theta)=(-1)^{s_\Gamma}C\Gamma_\mu (\theta, \eta)$, where $s_{\Gamma}=0$
if $C\Gamma_\mu$ is symmetric and $s_{\Gamma}=-1$ if $C\Gamma_\mu$ is skew-symmetric.
From these one can also find that
\be
C\Gamma_{\mu_1\dots\mu_p}(\eta, \theta)=(-1)^{{1\over2} p(p-1)} (-1)^{(p+1)s_C+p s_\Gamma}
C\Gamma_{\mu_1\dots\mu_p}(\theta, \eta)~.
\la{signpg}
\ee
Similarly we define
$
(\Gamma_{\mu_1\dots\mu_p}C^{-1})(\psi, \chi)= C^{-1}(\Gamma_{\mu_1\dots\mu_p}\psi,\chi)
$,
where $\chi,\psi\in \Delta^*$.

The product of two co-spinor  representations can be decomposed in terms of forms as
$\Delta^*\otimes \Delta^*=\sum^{n}_{p=1} \Lambda^p(V)\otimes \bC$.
In particular, one can write
\bea
(\psi\otimes \chi)(\eta\otimes\theta)
&=&{1\over {\rm dim}\Delta_n} \biggl( C^{-1}(\psi,\chi) C(\eta,\theta)
\cr
&+&\sum^{n}_{p=1} {(-1)^{p (s_\Gamma+s_C)}\over p!} (\Gamma^{\mu_1\dots\mu_p}C^{-1})(\psi,\chi)
~ C\Gamma_{\mu_1\dots\mu_p}(\eta,\theta)\biggr)~,
 \la{prod}
\eea
where $\eta,\theta\in \Delta$ and $\psi,\chi\in \Delta^*$.
The above decomposition is valid after restricting to $\Delta_\pm$ co-spinor
representations and to real co-spinor representations. We shall state the formulae later.
The formula (\ref{prod}) is also known as Fierz identity.

In the above formalism, it is possible to explicitly present the
spinors that are invariant under the action of certain subgroups
of $Spin(n)$. We shall mainly focus on the subgroups $G\subset Spin(n)$ which arise
as special holonomy groups in the Berger classification and the associated
manifolds admit a parallel spinor. These spinors have been given in \cite{ wang}. Here we shall summarize
the results and adjust the various formulae because of differences in the conventions.

(i) $G=SU(m)\subset Spin(2m)$. The invariant spinors under the
 $SU(m)\subset Spin(2n)$ are
 $1, e_1\wedge e_2\wedge\dots\wedge e_m$. This can be easily seen
 by decomposing the $Spin(2m)$ representations $\Delta^\pm$ under $SU(m)$.
 If $m=4k, 4k+2$ both invariant spinors are of the same chirality,
 ie they are elements of $\Delta^+$ while if $m=4k+1, 4k+3$, they have
 opposite chiralities. In addition observe that $\Gamma_j-i\Gamma_{m+j}(1)=0$
 and $\Gamma_j+i\Gamma_{m+j}(e_1\wedge\dots e_m)=0$, $j=1,\dots,m$. Therefore, the invariant
 spinors are pure spinors with respect to  the holomorphic and antiholomorphic parts of
 the decomposition of $V\otimes\bC$ with respect
 to the complex structure $J$.

 (ii) $G=Sp(k)\subset Spin(4k)$. The invariant spinors are
 $1,e_1\wedge e_2\wedge\dots\wedge e_{2k}, \omega, \omega^2,\dots, \omega^{k-1}$
 where $\omega=e_1\wedge e_2+\dots+e_{2k-1}\wedge e_{2k}$ which
 is the symplectic form in $U\subset \Delta^+$. Therefore there
 are $k+1$ parallel spinors.

 (iii) $G=Spin(7)\subset Spin(8)$. The invariant spinor is
 ${1\over \sqrt{2}} (e_1-e_2\wedge e_3\wedge e_4)$.

 (iv) $G=G_2\subset Spin(7)$.  The invariant spinor is
 ${1\over \sqrt{2}} (e_1-e_2\wedge e_3\wedge e_4)$.

\newsection{Differential and algebraic operations on spinors}

\subsection{First order differential operators}

Let $M$ be a spin manifold equipped with a spin connection $\nabla$
which admits a parallel spinor $\zeta$, $\nabla\zeta=0$. We shall focus on even-dimensional manifolds. Some
of the results can be easily extended to the odd dimensional case.
 We define $\eu_\pm=\Lambda^*(\Delta_\pm)$ and
$\eu=\Lambda^*(\Delta^*)$ equipped with the wedge product $\bar\wedge$.

\begin{definition}
The spin operator $d$ is a linear differential operator $d: \eu(M)\rightarrow \eu(M)$, and similarly
$d:\eu_\pm(M)\rightarrow \eu_\pm(M)$, such that
\be
d\phi=\sum_{\mu=1}^n C_\zeta^\mu \bar\wedge
\nabla_\mu \phi~,
\la{spinop}
\ee
where $C_\zeta^\mu(\eta)=C\Gamma^\mu(\zeta, \eta )$.
\end{definition}

Clearly $d: \eu^\ell\rightarrow \eu^{\ell+1}$ and
 $d: \eu_\pm^\ell\rightarrow \eu_\pm^{\ell+1}$. Choosing a basis in the
 space of co-spinors $\{\e^A: A=1, \dots, 2^m\}$, the $d$ operator
 can be written as
 \bea
 d\phi_{A_1 A_2\dots A_{\ell+1}}=\zeta^B C\Gamma^\mu_{BA_1}
 \nabla_\mu \phi_{A_2\dots A_{\ell+1}}&+&
 {\rm cyclic} (A_1, A_2,\dots, A_{\ell+1})
 \cr
 &&
 B, A_1,\dots, A_{\ell+1}=1, \dots, {\rm dim}\Delta~.
 \eea
 The operator $d$ depends on the choice of parallel spinor $\zeta$ and the connection
 $\nabla$. Although subject to the data above, the operator $d$ can always be defined on $\eu(M)$, the restriction
 on $d$ onto $\eu_+$ or $\eu_-$ depends on the choice of the parallel spinor $\zeta$.
 Since this depends on the dimension of the manifold and the choice
  of the parallel spinor, we shall explain the general properties of the operator $d$
  acting on $\eu$, and later we shall specialize into the various cases.

Evaluating $d^2$, we find
\be
d^2\phi={1\over2} C_\zeta^\mu \bar\wedge C_\zeta^\nu\bar\wedge
R_{\mu\nu} \phi~,
\ee
where $R$ is the curvature of the connection $\nabla$. Under certain conditions
the operator $d$ can be nilpotent, $d^2=0$. This depend on the choice of the spinor
$\zeta$ and the connection $\nabla$. There are two large classes of examples for which $d^2=0$.

\begin{itemize}

\item Group manifolds equipped with the left or the right invariant connections.

\item Manifolds that admit a pure parallel spinor.

\end{itemize}

In the case of group manifolds $R=0$. Therefore  $d$ operators
associated with the left- or right-invariant connections are all nilpotent.
 The different $d$ operators that
can be constructed  on group manifolds are determined by the orbits of $Spin(n)$ in $\Delta$.

 A spinor is pure if the subspace
\be
W(\zeta)=\{ v\in V_{\bC}~, v_\mu \Gamma^\mu\zeta=0\}
\ee
of $V_{\bC}$
has dimension ${1\over2} {\rm dim}_{\bC}(V_{\bC})$. We can use the inner product to decompose  $V_{\bC}=W(\zeta)\oplus Z$.

\begin{proposition}
The operator $d$ is nilpotent if, in addition to $\nabla\zeta=0$, the curvature
$R$ vanishes along the subspace
$\Lambda^2(Z)\subset \Lambda^2(V_{\bC})$, ie
\be
R|_{\Lambda^2(Z)}=0~. \la{nspin}
\ee
\end{proposition}
\bproof
The curvature $R$ can be viewed as a map from $R:\Lambda^2(M)\rightarrow \Lambda^2(M)$. Therefore
\be
d^2\phi={1\over2} C_\zeta^\mu\bar\wedge C_\zeta^\nu\bar\wedge  R_{\mu\nu}\phi=
{1\over2} [C_\zeta^\mu\bar\wedge C_\zeta^\nu\bar\wedge  R_{\mu\nu}]|_{\Lambda^2(Z)}\phi=0~.
\ee
\eproof

Clearly, this condition can be generalized to spinors $\zeta$ which are not pure but $W(\zeta)\not=\emptyset$.

The conditions on the curvature required for $d^2=0$ can also be determined
using (\ref{prod}). In particular we have

\begin{proposition}
The conditions on the curvature $R$ for $d^2=0$ can be expressed in terms of the
forms associated with the parallel spinor $\zeta$.
\end{proposition} \la{thone}
\bproof
We compute $d^2$ using (\ref{prod}) to find
\bea
d^2\phi&=& {1\over 2{\rm dim}\Delta_n} \biggl(
\sum^{n}_{p=0} {(-1)^{p (s_\Gamma+s_C)}\over p!} (\Gamma^{\rho_1\dots\rho_p}C^{-1})
(C_\zeta^\mu,C_\zeta^\nu)\biggr)
~ C\Gamma_{\rho_1\dots\rho_p} \bar\wedge R_{\mu\nu}\phi
\cr
&=&{1\over 2{\rm dim}\Delta_n} \biggl( \sum^{n}_{p=0}
{(-1)^{(p+1) (s_\Gamma+s_C)}\over p!} C(\zeta,
\Gamma^\mu\Gamma^{\rho_1\dots\rho_p}\Gamma^\nu\zeta)\biggr)
C\Gamma_{\rho_1\dots\rho_p}\bar\wedge R_{\mu\nu}\phi
\cr
&=&{1\over 2{\rm dim}\Delta_n} \biggl(\sum^{n}_{p=0} {(-1)^{(p+1) (s_\Gamma+s_C)}\over p!}
[C(\zeta,
\Gamma^{\mu\rho_1\dots\rho_p\nu}\zeta)
+ p (p-1)
g^{\mu\rho_1} C(\zeta,
\Gamma^{\rho_2\dots\rho_{p-1}}\zeta)  g^{\rho_p\nu}]\biggr)
\cr &&~~~~ C\Gamma_{\rho_1\dots\rho_p}
\bar\wedge R_{\mu\nu}\phi~,
\la{curvcon}
\eea
where $g$ is the metric on the manifold. In the above sum over $p$ only the
terms with $C\Gamma_{\rho_1\dots\rho_p}$ skew-symmetric contribute. Sufficient
conditions on the curvature for $d^2=0$ are
\be
[C(\zeta,
\Gamma^{\mu\rho_1\dots\rho_p\nu}\zeta)
+ p (p-1)
 C(\zeta,
\Gamma^{\rho_2\dots\rho_{p-1}}\zeta) g^{\mu\rho_1} g^{\nu\rho_p}] R_{\mu\nu}=0~
\la{curvconb}
\ee
for ${1\over2} p (p-1)+ (p+1)s_C+p s_\Gamma\in 2\bZ+1$.  In some cases these
conditions are also necessary.
\eproof

Provided that the condition for $d^2=0$ are met, we can define a cohomology theory associated
with the linear differential operator $d$.

\begin{definition}
The spin cohomology, $H_{d}(\eu)$, is that of the graded complex $(\eu, d)$, $d^2=0$, where $d$ is as in (\ref{spinop}).
Similarly, the spin cohomology, $H_d(\eu_\pm)$, is that of the graded complex $(\eu_\pm,d)$.
\end{definition}

\subsection{Twisted Complexes}

There are several ways to twist the complexes $\eu$ and $\eu_\pm$.
Here we shall consider two cases which we shall describe below.

\subsubsection{The complexes $\eu\otimes E$ and $\eu_\pm\otimes E$}

Let $E$ be a vector bundle $E$ over the spin manifold $M$ equipped
with a connection $\nabla^E$. One way to  twist the complexes $\eu$ and $\eu_\pm$
is  to consider $\eu\otimes E$ and $\eu_\pm\otimes E$.  Let $\zeta$ be a parallel
spinor with respect to a spin connection $\nabla^M$ on the manifold $M$
induced from the tangent bundle, $\nabla^M\zeta=0$. The spin differential operator $d$ is
\be
d\phi=C^\mu_\zeta\bar\wedge \nabla_\mu\phi~,~~~~~~
\la{twtwe}
\ee
where $\nabla=\nabla^M\otimes 1+1\otimes \nabla^E$ on $\eu\otimes E$ or on $\eu_\pm\otimes E$
and $\phi\in \eu\otimes E$ or $\eu_\pm\otimes E$, respectively.

The condition $d^2=0$ implies conditions on both the curvature $R$ of $M$
and the curvature $F$ of the connection $\nabla^E$ of the bundle $E$.

\begin{theorem}
The operator $d$ is nilpotent providing that both the curvature $R$
of the manifold $M$ and the curvature $R$ of $E$ satisfy either (\ref{nspin}) or (\ref{curvcon}).
\end{theorem}
\bproof
This is similar to the proof given in the previous section.
\eproof

There is a particular twisted complex of the this type that we shall consider
by taking $E=\Lambda^*(M)$ or $E=\Lambda^*(M)\otimes \bC$. We shall see
that in this case one can define certain algebraic operators with are nilpotent.
The spin cohomology, $H_d(\eu\otimes E)$, of the linear operator $d$ for the twisted complex $(\eu\otimes E, d)$,
 can be defined in analogy
with the spin cohomology of the untwisted case in the previous section. This definition can be extended
for $H_d(\eu_\pm\otimes E)$.

\subsubsection{The complexes $\eu(E)=\Lambda^*(\Delta^*\otimes E)$ and
$\eu(E)_\pm=\Lambda^* (\Delta_\pm\otimes E)$}

These complexes allow the definition of the spin operator $d$ operator on manifolds
that do not admit a spin structure but admit a $Spin_c$ or in general a
$Spin_G$ structure.
Another use of twisted complexes  $\eu(E)$ is that they allow the imposition of a reality condition.
It is known that there are not  real (Majorana) spin
representations for $n=8k+4$ dimensional manifolds and so there is not a real
complex $\eu$. However it is
possible to construct a real complex $\eu(E)$ by taking $E$ to be a rank two $SU(2)$ bundle.

Let $\zeta$ be a parallel section of $\Delta\otimes E$ with respect to a connection
${\cal D}=\nabla^M\otimes 1+1\otimes \nabla^E$, where $\nabla^M$ is a spin connection on the manifold $M$
induced from the tangent bundle and $\nabla^E$ is a connection on the vector bundle $E$.
The operator $d$ on $\eu(E)$ or $\eu(E)_\pm$ is defined as
\be
d\phi=C^\mu_\zeta\bar\wedge {\cal D}_\mu\phi~,~~~~~~
\ee
where $\bar\wedge$ is the wedge operation in $\Lambda^*(\Delta^*\otimes E)$.
One consequence of this definition is that $\zeta$ is not necessary a parallel
section of the spinor bundle but of $\Delta\otimes E$.

The case that it is of most
interest to us is that for which $E$ is a line bundle. In this case, the conditions
for $d^2=0$ can be expressed as conditions on the curvature $R$ and $F$ of the manifold
and of the line bundle, respectively. The formulae are similar to those
in (\ref{nspin}) and (\ref{curvcon}).

The construction can be further generalized in the case for which there is no a spin structure
but there is a $Spin_c$ structure. In this case although the spin bundle $\Delta$ is not well-defined
$\Delta\otimes E$ is and so is $\eu(E)$.

Provided that the conditions for $d^2=0$ are met, we can define the twisted spin cohomology, $H_d(\eu(E))$, of
the graded complex $(\eu(E), d)$. Similary we can define the twisted spin cohomology, $H_d(\eu_\pm(E))$, of
the graded complex $(\eu_\pm(E), d)$

\subsection{Algebraic operations}

\subsubsection{The algebraic operator $D_{(p)}$}

There are several  algebraic cohomology operations that can be defined
on the twisted complexes $\Lambda^*\otimes \eu$, $\Lambda^*\otimes \eu_\pm$,
$Sym^*\otimes \eu$ and $Sym^*\otimes \eu_\pm$, where $Sym^*=\oplus^\infty_{p=0} Sym^p$
and $Sym^p$ is the symmetrized product of $p$ copies of $\Lambda^1$.

The maps $C\Gamma^{(p)}: \Delta\otimes \Delta\rightarrow \Lambda^p$
\be
C\Gamma^{(p)}(\eta,\theta)={1\over p!}C\Gamma_{\mu_1\dots\mu_p}(\eta, \theta) e^{\mu_1}\wedge\dots\wedge e^{\mu_p}
\ee
are skew-symmetric, ie
$C\Gamma^{(p)}(\eta,\theta)=- C\Gamma^{(p)}(\theta,\eta)$, provided
that ${1\over 2} p (p-1)+(p+1) s_C+p s_\Gamma\in 2\bZ+1$ as it can been seen from (\ref{signpg}).

\begin{definition}
 The algebraic spin operator
$D_{(p)}: \Lambda^q(M)\otimes \eu^\ell(M)\rightarrow \Lambda^{q-p}(M)
\otimes \eu^{\ell+2}(M)$
is
\bea
D_{(p)}\phi={(-1)^{{1\over2} p (p-1)+\ell}\over2~ (q-p)! p! \ell!} (C\Gamma^{\mu_1\dots\mu_p})_{A_1 A_2}
\phi_{\mu_1\dots\mu_p \nu_1\dots \nu_{q-p} A_3\dots A_{\ell+2}}
\cr
 e^{\nu_1}\wedge\dots
\wedge e^{\nu_{q-p}}\otimes \epsilon^{A_1}\bar\wedge\dots\bar\wedge \epsilon^{A_{\ell+2}},
\eea
if $q\geq p$ and $D_{(p)}=0$ for $p>q$, where $C\Gamma^{(p)}$ is skew-symmetric.
\end{definition}

It is straightforward to show that

\begin{proposition}
$D_{(p)}$ is nilpotent, $D_{(p)}^2=0$, provided that $p\in 2\bZ+1$.
\end{proposition}

\begin{proposition}
The algebraic spin operator $D_{(p)}$ can be restricted on $\Lambda^*\otimes \eu_\pm$, iff ~${\rm dim}~M=8k+2, 8k+6$.
\end{proposition}
\bproof
It can be seen from the properties of spinor inner product $C$ summarized in section two that for ${\rm dim}~M=8k+2, 8k+6$,
$C\Gamma^{(p)}: \Delta_\pm\otimes \Delta_\pm\rightarrow \Lambda^p$, $p\in 2\bZ+1$.
\eproof

In what follows, when we refer to the algebraic spin operator on $\Lambda^*\otimes \eu_\pm$ complexes we shall assume
the condition of the above proposition applies and ${\rm dim}~M=8k+2, 8k+6$.

The algebraic spin operator $D_{(p)}$ can be extended to twisted complexes
$\Lambda^*\otimes \eu \otimes E$ and $\Lambda^*\otimes \eu_\pm \otimes E$ in a straightforward way.
There is also an extension to the twisted complexes $\Lambda^*\otimes \eu(E)$ and $\Lambda^*\otimes\eu_\pm(E)$ provided
that $E$ is equipped with and inner product $h$. In particular we define
\bea
D_{(p)}:&& \Lambda^q(M)\otimes \eu^\ell(E)\rightarrow \Lambda^{q-p}(M)
\otimes \eu^{\ell+2}(E)
\eea
where
\bea
D_{(p)}\phi={(-1)^{{1\over2} p (p-1)+\ell}\over2~ (q-p)! p! \ell!}
(C\Gamma^{\mu_1\dots\mu_p}\otimes h)_{I_1A_1, I_2A_2}
\phi_{\mu_1\dots\mu_p \nu_1\dots \nu_{q-p} I_3A_3,\dots, I_{\ell+2}A_{\ell+2}}
\cr
 e^{\nu_1}\wedge\dots
\wedge e^{\nu_{q-p}}\otimes \epsilon^{I_1A_1}\bar\wedge\dots\bar\wedge
\epsilon^{I_{\ell+2}A_{\ell+2}},
\la{talg}
\eea
if $q\geq p$ and $D_{(p)}=0$ for $p>q$. This operation is well defined
provided  that $C\Gamma^{(p)}\otimes h$ is skew-symmetric.
This is the case when either  $C\Gamma^{(p)}$ is symmetric and
$h$ is skew-symmetric or $C\Gamma^{(p)}$ is skew-symmetric and
$h$ is symmetric. In all the above cases $D_{(p)}^2=0$ provided $p\in 2\bZ+1$.

\begin{definition}
The algebraic spin cohomology $H_{D_{(p)}}(\Lambda^*\otimes \eu)$ is defined as the cohomology of the
double graded complex $(\Lambda^*\otimes \eu, D_{(p)})$. This definition can be extended to the
rest of the twisted and untwisted spin complexes.
\end{definition}

A particular case of this operation is for $p=1$. In this case, we have
$D=D_{(1)}: \Lambda^*\otimes \eu\rightarrow \Lambda^*\otimes\eu$,
where
\bea
D\phi={(-1)^{\ell}\over2~ (q-1)!  \ell!} (C\Gamma^{\mu})_{A_1 A_2}
\phi_{\mu_1 \nu_1\dots \nu_{q-1} A_3\dots A_{\ell+2}}
 e^{\nu_1}\wedge\dots
\wedge e^{\nu_{q-1}}\otimes \epsilon^{A_1}\bar\wedge\dots\bar\wedge \epsilon^{A_{\ell+2}},
\eea
if $q\geq 1$ and $D_{(p)}=0$ for $q=0$.
In particular $D$ is defined on $\Lambda^*\otimes \eu$ and $\Lambda^*\otimes \eu\otimes E$
for $m=4k, 4k+2$ if $C=A$ and for $m=4k+3$ if $C=B$. It is also defined on $\Lambda^*\otimes \eu_\pm$ and
$\Lambda^*\otimes \eu_\pm\otimes E$ for $m=4k+2$ if $C=A$ and for $m=4k+3$ if $C=B$,~ ${\rm dim}~M=2m$.
In the twisted case $\Lambda^*\otimes \eu(E)$, the operator $D$ can be defined
in all the cases for which $E$ admits a fibre inner product such that
$C\Gamma^{(1)}\otimes h$ is skew-symmetric.

\begin{proposition}
$D_{(p)}$  anti-commutes with the differential operator $d$, ie
\bea
dD_{(p)}+D_{(p)} d&=&0~.
\eea
\end{proposition}
\bproof
We can show this after a direct computation using the property of the connection $\nabla^M$ of the manifold to be
a spin connection
induced from the tangent bundle. In the twisted case (\ref{talg}) this also the case provided that $\nabla^E h=0$.
\eproof

The differential and algebraic spin cohomology operators on the various untwisted and twisted complexes above
can be combined into an new spin cohomology
operator $d+D$. The new cohomology  operator $d+D$ defines a new cohomology,  $H_{d+D}$, which can be computed
using  spectral sequences. We shall describe such  computation on
Calabi-Yau manifolds of dimension six.

\subsubsection{The algebraic operator $\hat D$}

Apart from the $D_{(p)}$ algebraic operator, there is another algebraic operation  $\hat D$.

\begin{definition}
The algebraic operator $\hat D$  on the complex
$Sym^*\otimes \eu$ is
\bea
\hat D:&& Sym^q(M)\otimes \eu^\ell\rightarrow Sym^{q-1}(M)\otimes \eu^{\ell+1}~,
\eea
where
\bea
 \hat D\phi=(-1)^\ell{1\over (q-1)! \ell!}(C_\zeta)_{A_1}^\mu  \phi_{\mu\nu_1\dots\nu_{q-1} A_2\dots A_{\ell+1}}
 e^{\nu_1}\wedge\dots
\wedge e^{\nu_{q-1}}\otimes \epsilon^{A_1}\bar\wedge\dots\bar\wedge \epsilon^{A_{\ell+1}}~
\la{algtwo}
\eea
if $q\geq 1$ and $\hat D=0$ if $q=0$.
\end{definition}

It is straightforward to extend this definition to the other untwisted and twisted complexes.
Moreover one can show that
\begin{proposition}
$\hat D$ is nilpotent, $\hat D^2=0$.
\end{proposition}

As in the case of $D_{(p)}$ algebraic spin operator
\begin{proposition}
 $\hat D$ anti-commutes with the differential operators $d$, ie
\bea
d\hat D+\hat D d&=&0~.
\eea
\end{proposition}

A consequence of this is that one can define a new cohomology operator $d+\hat D$ and an associated
cohomology $H_{d+\hat D}$ which can be computed using a spectral sequence.

\begin{proposition}
Let $\hat D$ be the spin algebraic operator on $Sym^*\otimes \eu$ defined as in (\ref{algtwo}).
If $C_\zeta$ is an isomorphism, then $H_{\hat D}^*=\eu$.
\end{proposition}
\bproof
Since $C_\zeta$ is an isomorphism, then $Sym^*\otimes \eu=Sym^*\otimes \Lambda^*$.
$Sym^p\otimes \Lambda^q$ can be decomposed under $GL(n, \bC)$ into two irreducible
representations. These have dimensions
\bea
\delta_1&=&{n (n+1)\dots (n+p-1) (n-1)\dots (n-q)\over (p+q) (p-1)! q!}
\cr
\delta_2&=&{n (n+1)\dots (n+p) (n-1)\dots (n-q+1)\over (p+q) p! (q-1)!}~.
\eea
Clearly ${\rm Ker}\hat D|_{Sym^0\times \Lambda^q}=\Lambda^q$. In addition
\be
{\rm Ker}\hat D|_{Sym^p\times \Lambda^q}=\hat D({Sym^{p+1}\times \Lambda^{q-1}}),~~~~p>0~,
\ee
with ${\rm dim}~{\rm Ker}\hat D|_{Sym^p\times \Lambda^q}=\delta_2$
Therefore all cohomology $H^{p,q}_{\hat D}=0$ for $p>0$ and $H^{0,q}_{\hat D}=\Lambda^q=\eu^q$.
\eproof

This theorem can be thought as a consequence of the Spencer cohomology \cite{spencer}.
Clearly the above result can be generalized to $Sym^*\otimes \eu_\pm$ and twisted complexes.

\newsection{Manifolds with connections of holonomy  $SU(m)$ and spin cohomology}

As we have mentioned on even-dimensional
Riemannian manifolds, there are two complex spin representations $\Delta^\pm$. In addition
there are real spin representations provided that $m=4k+1, 4k+3, 4k$, ${\rm dim}~M=2m$. In what follows,
we shall focus on the spin cohomology associated  the complex representations. The spin cohomology
associated with real representations will be investigated later.

\subsection{Complex spinor representations}

Let $M$ be a Riemannian manifold equipped with a spin connection
$\nabla$ with ${\rm hol}(\nabla)\subseteq SU(m)$.
We take that the metric on $M$ to be compatible with
the parallel almost complex structure $J$.
There are two distinct $\nabla$-parallel complex spinors. These are given by
$\zeta_1=1, \zeta_2=e_1\wedge\dots\wedge e_m$. These spinors are
 of different chirality if $m=4k+1, 4k+3$ and of the
same chirality if $m=4k, 4k+2$. Therefore there are two first order differential spin
operators $d_1$ and $d_2$ associated with the spinors $\zeta_1$ and $\zeta_2$, respectively.
If $m$ is odd, $d_1:\eu_+\rightarrow \eu_+$ and $d_2:\eu_-\rightarrow \eu_-$
while if $m$ is even, $d_1, d_2: \eu_-\rightarrow \eu_-$. We shall treat the two
cases separately.

\subsubsection{The $m=4k+1, 4k+3$ case}
\begin{theorem}
The operator $d_1: \eu_+\rightarrow \eu_+$ is nilpotent, $d_1^2=0$,
provided that the (2,0) part
of the curvature $R$ of the connection $\nabla$ with respect to
the almost complex structure $J$ vanishes.
\end{theorem}
\bproof
For this we compute $d_1^2$ on $\eu_+^\ell$ to find
\be
d_1^2 \phi={1\over2} C_1^\mu\bar\wedge C_1^\nu\bar\wedge R_{\mu\nu} \phi~,
\ee
where  $R_{\mu\nu}=[\nabla_\mu, \nabla_\nu]$ is the curvature of the connection $\nabla$.
In a spinor basis $\{\epsilon^{a}: a=1,\dots,{\rm dim} \Delta_+$,
the above expression can be written as
\be
d_1^2\phi={1\over2~~\ell!}
(C_1)^\mu_{a_1} (C_1)_{a_2}^\nu (R_{\mu\nu} \phi)_{a_3\dots a_{\ell+2}}~\epsilon^{a_1}
\wedge\epsilon^{a_2}\wedge\epsilon^{a_3}\wedge\dots \wedge \epsilon^{\ell+2}~.
\ee
Observe that the product representation $\Delta_+\otimes \Delta_+$
can be decomposed as
\be
\Delta_+\otimes \Delta_+=
\sum^{{m-3\over2}}_{p=1} \Lambda^{2p+1}(V_{\bC})\oplus \Lambda^{m+}(V_{\bC})~.
\ee
In particular using (\ref{prod}), we have
\be
\chi\otimes \psi(\eta \otimes \theta)=
{1\over {\rm dim}_{\bC}\Delta_+}~\sum^{{m-3\over2}}_{p=1}
 {(-1)^{(2p+1)(s_\Gamma+s_C)}\over (2p+1)!}
(\Gamma^{\mu_1\dots \mu_{2p+1}}C^{-1})(\chi,\psi)
C\Gamma_{\mu_1\dots \mu_{2p+1}}(\eta,\theta)~.
\la{proda}
\ee
The only non-vanishing form associated with  $\zeta_1$ spinor is
the m-form given by
\be
\epsilon_1=  {1\over m!}
C(\zeta_1, \Gamma_{\rho_1\dots\rho_m}\zeta_1) e^{\rho_1}\wedge\dots\wedge e^{\rho_m}
~.
\ee
Moreover, from section two, we have that $C\Gamma_{\rho_1\dots\rho_m}$ is
symmetric and $C\Gamma_{\rho_1\dots\rho_{m-2}}$
is skew-symmetric. Applying the formula (\ref{proda})
for $\chi=C_1^\mu$ and $\psi=C^\nu_1$, we find that
the only non-vanishing term is
\bea
d_1^2 \phi&=&{1\over2} {(-1)^{(m-2)(s_\Gamma+s_C)}\over (m-2)! {\rm dim}_{\bC}\Delta_+}
(\Gamma^{\rho_1\dots \rho_{m-2}}C^{-1})(C_1^\mu, C_1^\nu) C\Gamma_{\rho_1\dots\rho_{m-2}}\bar\wedge
R_{\mu\nu}\phi
\cr
&=&
{1\over2} {(-1)^{(m-1)(s_\Gamma+s_C)}\over (m-2)! {\rm dim}_{\bC}\Delta_+}
C(1, \Gamma^{\mu\rho_1\dots\rho_{m-2}\nu}1)
C\Gamma_{\rho_1\dots\rho_{m-2}}\bar\wedge R_{\mu\nu}\phi~.
\la{dab}
\eea
Therefore $d_1^2=0$, if  the (2,0) component of the curvature  $R$ vanishes since
$(\Gamma_j-i\Gamma_{j+m})1=0$ and so $\epsilon_1$ is an $(0,m)$ form.
\eproof

\begin{theorem}
The operator $d_2: \eu_-\rightarrow \eu_-$ is nilpotent, $d_2^2=0$,
provided that the (0,2) part
of the curvature $R$ of the connection $\nabla$
with respect to the almost complex structure $J$ vanishes.
\end{theorem}
\bproof
The proof of this is similar to the one presented above for the case of the
$\zeta_1$ parallel spinor. One difference is the decomposition
\be
\Delta_-\otimes \Delta_-
=\sum^{{m-3\over2}}_{p=1} \Lambda^{2p+1}(V_{\bC})\oplus \Lambda^{m-}(V_{\bC})~.
\ee
A direct computation reveals that
\bea
d_2^2 \phi&=&{1\over2} {(-1)^{(m-2)(s_\Gamma+s_C)}\over (m-2)! {\rm dim}_{\bC}\Delta_+}
(\Gamma^{\rho_1\dots \rho_{m-2}}C^{-1})(C_2^\mu, C_2^\nu) C\Gamma_{\rho_1\dots\rho_{m-2}}\bar\wedge
R_{\mu\nu}\phi
\cr
&=&
{1\over2} {(-1)^{(m-1)(s_\Gamma+s_C)}\over (m-2)! {\rm dim}_{\bC}\Delta_+}C(e_1\wedge\dots\wedge
e_m, \Gamma^{\mu\rho_1\dots\rho_{m-2}\nu}e_1\wedge\dots\wedge e_m)
\cr
~~~~~~~~~~~~&&
C\Gamma_{\rho_1\dots\rho_{m-2}}\bar\wedge R_{\mu\nu}\phi~.
\la{dbb}
\eea
Using $(\Gamma_j+i\Gamma_{j+m})e_1\wedge\dots\wedge e_m=0$,  we conclude
that $d_2^2=0$ if the (0,2) part of the curvature $R$ vanishes.
\eproof

\begin{corollary}
The operator $d=d_1\oplus d_2: \eu_+\oplus \eu_-\rightarrow \eu_+\oplus \eu_-$ is
nilpotent provided that
the curvature $R$ of the connection $\nabla$ is (1,1)
with respect to the almost complex
structure $J$.
\end{corollary}
\bproof
It follows immediately from  the two theorems above.
\eproof

We therefore conclude that there are three kinds of
untwisted differential spin cohomology associated
with a manifold that admits a connection with holonomy contained in $SU(m)$,
$m=4k+1, 4k+3$. The complexes are $(\eu_+, d_1)$, $(\eu_-, d_2)$ and $(\eu_+\oplus
\eu_-, d_1\oplus d_2)$ and the associated spin cohomologies are denoted  as
$H_{d_1}(\eu_+)$, $H_{d_2}(\eu_-)$ and $H_d(\eu_+\oplus\eu_-)$,
respectively.

\subsubsection{The $m=4k, 4k+2$ case}

In this case  both parallel spinors $\zeta_1, \zeta_2\in \Delta^+$. Therefore
 $d_1, d_2: \eu_-\rightarrow \eu_-$.
\begin{theorem}
The operators $d_1, d_2: \eu_-\rightarrow \eu_-$ are nilpotent, $d_1^2=0$ and $d_2^2=0$,
provided that the  either (2,0) or the (0,2) part
of the curvature $R$ of the connection $\nabla$ with respect to the
almost complex structure $J$ vanishes, respectively.
\end{theorem}
\bproof
The proof of this statement is similar to that  for the cases
$m=4k+1, 4k+3$ described in the previous section. In particular, we have
\be
\Delta_-\otimes \Delta_-
=\sum^{{m-2\over2}}_{p=0} \Lambda^{2p}(V_{\bC})\oplus \Lambda^{m-}(V_{\bC})~.
\la{dec}
\ee
{}From the results of section two, the map $C\Gamma^{\mu_1\dots\mu_m}$
is symmetric and $C\Gamma^{\mu_1\dots\mu_{m-2}}$ is skew-symmetric
with respect to both inner products $C=A,B$. The expressions for $d_1^2$ and $d_2^2$
are given by (\ref{dab}) and (\ref{dbb}), respectively. From these,
it is straightforward to see that $d_1^2=0$ ($d_2^2=0$) if the (0,2) ((2,0))
part of the curvature $R$ of $\nabla$ vanishes.
\eproof

Since both operators $d_1, d_2$ act on the same complex, one can define the
operator $d=d_1+d_2: \eu_-\rightarrow \eu_-$. If $d_1^2=d_2^2=0$, $d^2=d_1 d_2+d_2 d_1$.
Therefore $d$ is nilpotent iff the operators $d_1, d_2$ anti-commute.

\begin{theorem}
The operator $d_1 d_2+d_2 d_1= 0$  iff the curvature of $\nabla$ vanishes $R=0$.
\end{theorem}
\bproof
Applying the definition of $d_1$ and $d_2$, one can  find that
\bea
(d_1 d_2+d_2 d_1)\phi
= C_1^\mu\bar\wedge C_2^\nu \bar\wedge \nabla_\mu\nabla_\nu \phi
- C_1^\mu\bar\wedge C_2^\nu \bar\wedge  \nabla_\nu\nabla_\mu \phi
=C_2^\mu\bar\wedge C_1^\nu \bar\wedge R_{\mu\nu}\phi~.
\eea

The K\"ahler form associated with the parallel spinors is
\be
\Omega={-i\over2 C(\zeta_2, \zeta_1)} C(\zeta_2, \Gamma_{\mu\nu}\zeta_1)
e^\mu\wedge e^\nu~.
\ee
It can then be seen that
\be
C\Gamma_{(2p)}(\zeta_2, \zeta_1)={1\over (2p)!}C(\zeta_2, \Gamma_{\rho_1\dots\rho_{2p}}\zeta_1)
e^{\rho_1}\wedge\dots\wedge e^{\rho_{2p}}
={(-i)^p C(\zeta_2, \zeta_1)\over p!} \wedge^p\Omega~.
\ee
Applying (\ref{dec}), we find
\bea
(d_1 d_2+d_2 d_1)\phi&=&{1\over {\rm dim}_{\bC}
(\Delta_+)}\sum_{p=0}^{{m\over2}} {(-1)^{s_C+s_\Gamma}\over (2p)!}
C(\zeta_2, \Gamma^\mu \Gamma^{\rho_1\dots \rho_{2p}} \Gamma^\nu \zeta_1)
 C\Gamma_{\rho_1\dots \rho_{2p}}
\bar\wedge R_{\mu\nu}\phi
\cr
&=&{1\over {\rm dim}_{\bC}(\Delta_+)} \sum_{p=0}^{{m\over2}}
 {(-1)^{s_C+s_\Gamma}\over (2p)!}
\biggl[C(\zeta_2, \Gamma^{\mu\rho_1\dots \rho_{2p}\nu}  \zeta_1)
\cr
&+&(2p) (2p-1)g^{\mu \rho_1} C(\zeta_2, \Gamma^{\rho_2\dots \rho_{2p-1}}\zeta_1)
 g^{\rho_{2p}\nu}
\biggr]
 C\Gamma_{\rho_1\dots \rho_{2p}}
\bar\wedge R_{\mu\nu}\phi
~,
\eea
where $g$ is the metric on the manifold $M$.
This can be rewritten as
\begin{small}
\bea
(d_1 d_2+d_2 d_1)\phi&=&{C(\zeta_2, \zeta_1) \over {\rm dim}_{\bC}(\Delta_+)}
\sum_{p=0}^{{m\over2}} {(-1)^{s_C+s_\Gamma}\over (2p)!}
\cr
&&
\biggl[{(-i)^{p+1}
(2p+2)!\over 2^{p+1} (2p+1) (p+1)!} \Omega^{\mu\nu} \Omega^{\rho_1\rho_2}\dots
\Omega^{\rho_{2p-1}\rho_{2p}} C\Gamma_{\rho_1\dots \rho_{2p}}
\bar\wedge R_{\mu\nu}\phi
\cr
&+& {(-i)^{p-1} (2p)!\over 2^{p-1} (p-1)!} \Omega^{\rho_2\rho_3}\dots
\Omega^{\rho_{2p-2}\rho_{2p-1}}
\cr
&& C\Gamma_{\rho_1\dots \rho_{2p}}
\bar\wedge ( R_{\mu\nu}J^\mu{}_{\rho_1}J^\nu{}_{\rho_{2p}}+R_{\rho_1\rho_{2p}})\phi
\biggl]
\eea
\end{small}

{}For $m=4k$, $C\Gamma_{(2p)}$  is symmetric for $p=2q$ while
they are skew-symmetric
for $p=2q+1$. Therefore in this case only the latter terms contribute in the sum.
Similarly for $m=4k+2$, $C\Gamma_{(2p)}$  is skew-symmetric for $p=2q$
while they are symmetric
for $p=2q+1$. Therefore only the former terms contribute is the sum.

It is clear that if the (1,1) part of the curvature vanishes, then the proposition
is satisfied. However the (2,0) and the (0,2) parts of the curvature vanish as well.
Thus $R=0$.
\eproof

We therefore conclude that there are three kinds of
untwisted differential spin cohomology associated
with a manifold that admits a connection with holonomy contained in $SU(m)$,
$m=4k, 4k+2$. The complexes are $(\eu_-, d_1)$, $(\eu_-, d_2)$ and $(\eu_-, d=d_1+d_2)$
and the associated spin cohomologies are denoted  as
$H_{d_1}(\eu_-)$, $H_{d_2}(\eu_-)$ and $H_d(\eu_-)$,
respectively. Unlike the case case where $m=4k+1, 4k+3$, all three cohomologies
are cohomologies of the complex $\eu_-$.

\subsection{Adjoint operators and Laplacians}

As we have mentioned in section two, $\Delta_\pm$ are equipped with a $Spin(n)$-invariant
inner  product. Because of this, one can find the adjoints of the
 spin cohomology operators $d_1, d_2$ and their associated Laplacians.
 As in the previous section, we shall distinguish between the $m=4k+1$, $4k+3$ and
$m=4k$, $4k+2$ cases. This is because of the properties of the inner product
are different-see section two.

\subsection{The $m=4k+1$, $4k+3$ case}

We extend the inner product $C^{-1}$ from $\Delta_+\oplus\Delta_-$ to
the space of sections of $\eu_+\oplus \eu_-$
and denote it with the same symbol. The inner product $C^{-1}$  vanishes
if it is restricted on either $\eu_+$ or $\eu_-$.

\begin{definition}
The adjoint operator $\delta_1: \eu_-\rightarrow \eu_-$
of $d_1: \eu_+\rightarrow \eu_+$
is
\be
C^{-1}(\phi, d_1\psi)=C^{-1}(\delta_1\phi, \psi)~.
\ee
Similarly, the adjoint operator $\delta_2: \eu_+\rightarrow \eu_+$ of
$d_2: \eu_-\rightarrow \eu_-$
is
\be
C^{-1}(\psi, d_2\phi)=C^{-1}(\delta_2\psi, \phi)~.
\ee
\end{definition}

Using these adjoints,  one can define two Laplace operators
$\Delta_1=\delta_2 d_1+d_1 \delta_2$
and $\Delta_2=\delta_1 d_2+d_2 \delta_1$. The Laplace operator
of $d: \eu_+\oplus \eu_-\rightarrow
\eu_+\oplus \eu_-$ is $\Delta= \Delta_1\oplus \Delta_2$.

To compute the Laplace operator $\Delta_1$, we use  the above definitions to find
\bea
\Delta_1\phi=-(-1)^{s_C+s_\Gamma} \biggl(C(\zeta_1, \Gamma^\nu \Gamma^\mu\zeta_2)
\nabla_\mu \nabla_\nu \phi + (-1)^{s_c+s_\Gamma} C_1^\nu \bar\wedge
(C_2^\mu\rc R_{\nu\mu}\phi)\biggr)~,
\la{lapa}
\eea
where $\eta\rc \phi$ denotes inner derivation with respect to spinor $\eta$, ie
\be
\eta\rc \phi= {(-1)^{s_c}\over \ell!} \eta^B \phi_{BA_1\dots A_\ell}
 \epsilon^{A_1}\bar\wedge\dots\bar\wedge
\epsilon^{A_{\ell}}= {(-1)^{s_c}\over \ell!}( C^{-1})^{BE}\eta_B
\phi_{EA_1\dots A_\ell} \epsilon^{A_1}
\bar\wedge\dots\bar\wedge
\epsilon^{A_{\ell}}~,
\ee
which is equivalent to
\be
(\eta\rc \phi)_{A_1\dots A_{\ell}}=(\ell+1) (-1)^{s_c} \eta^B \phi_{BA_1\dots A_\ell}~.
\ee
More generally, we have
\be
(\eta\rc \phi)_{E_1\dots E_q,A_1\dots A_{\ell}}
=(\ell+1) (-1)^{s_c} \eta_{E_1\dots E_q}{} ^B \phi_{BA_1\dots A_\ell}
\ee

The product of the co-spinor representations $\Delta_\pm$ can be decomposed as
\be
\Delta_+\otimes \Delta_-=\sum_{p=0}^{[{m\over2}]} \Lambda^{2p}~.
\ee
The formula that relates the product of two co-spinors to forms is given by  (\ref{prod}) after
the appropriate restrictions.
Applying this to the second term in the Laplace operator, we find
\bea
\Delta_1\phi&=&-(-1)^{s_C+s_\Gamma} \biggl( C( \zeta_1, \Gamma^\nu \Gamma^\mu\zeta_2)
\nabla_\mu \nabla_\nu \phi
\cr
&+& {1\over {\rm dim}_{\bC}(\Delta_+)}\sum_{p=0}^{[{m\over2}]}
{1\over (2p)!} C(\zeta_1, \Gamma^\mu \Gamma^{\rho_1\dots \rho_{2p}} \Gamma^\nu \zeta_2)
 C\Gamma_{\rho_1\dots \rho_{2p}}
\rc R_{\mu\nu}\phi\biggr)~.
\eea
In turn this can be written as
\bea
\Delta_1\phi&=&-(-1)^{s_C+s_\Gamma} C(\zeta_1, \zeta_2) \biggl(
 (g^{\mu\nu}+i\Omega^{\mu\nu})
\nabla_\mu \nabla_\nu \phi
\cr
&+& {1\over{\rm dim}_{\bC}(\Delta_+)}
\sum_{p=0}^{[{m\over2}]} {1\over (2p)!}
\cr
&&
\biggl[{i^{p+1}
(2p+2)!\over 2^{p+1} (2p+1) (p+1)!} \Omega^{\mu\nu} \Omega^{\rho_1\rho_2}\dots
\Omega^{\rho_{2p-1}\rho_{2p}} C\Gamma_{\rho_1\dots \rho_{2p}}
\rc R_{\mu\nu}\phi
\cr
&+& {i^{p-1} (2p)!\over 2^{p-2} (p-1)!} \Omega^{\rho_2\rho_3}\dots
\Omega^{\rho_{2p-2}\rho_{2p-1}} g^{\rho_1 \mu} g^{\rho_{2p} \nu} C\Gamma_{\rho_1\dots \rho_{2p}}
\rc R^{1,1}_{\mu\nu}\phi
\biggl]\biggr)~,
\eea
where
\be
R^{1,1}_{\rho_1\rho_{2p}}={1\over2}(
 R_{\mu\nu} J^{\mu}{}_{\rho_1} J^{\nu}{}_{\rho_{2p}}+R_{\rho_1\rho_{2p}})~.
 \ee
 is the $(1,1)$ part of the curvature with respect to the almost complex structure $J$.

The Laplace operators $\Delta_2$ is given as in (\ref{lapa})
 but with the parallel spinors
$\zeta_1$ and $\zeta_2$ interchanged. The effect that this can be easily
computed from \ref{lapa}) using the symmetry properties of $C\Gamma_{\mu_1\dots\mu_q}$.
 In particular we find that

\begin{corollary}
The $\Delta_2$ Laplace operator is
\bea
\Delta_2\phi&=&-(-1)^{s_C+s_\Gamma} C(\zeta_1, \zeta_2) \biggl( (-1)^{s_c}
 (g^{\mu\nu}-i\Omega^{\mu\nu})
\nabla_\mu \nabla_\nu \phi
\cr
&+& {(-1)^{s_c}\over{\rm dim}_{\bC}(\Delta_+)}
\sum_{p=0}^{[{m\over2}]}
\biggl[{(-i)^{p+1}
(2p+2)!\over 2^{p+1} (2p+1) (p+1)!} \Omega^{\mu\nu} \Omega^{\rho_1\rho_2}\dots
\Omega^{\rho_{2p-1}\rho_{2p}} C\Gamma_{\rho_1\dots \rho_{2p}}
\rc R_{\mu\nu}\phi
\cr
&+& {(-i)^{p-1} (2p)!\over 2^{p-2} (p-1)!} \Omega^{\rho_2\rho_3}\dots
\Omega^{\rho_{2p-2}\rho_{2p-1}} g^{\rho_1 \mu} g^{\rho_{2p} \nu} C\Gamma_{\rho_1\dots \rho_{2p}}
\rc R^{1,1}_{\mu\nu}\phi
\biggl]\biggr)~.
\la{deltatwo}
\eea
\end{corollary}

\subsection{The $m=4k$, $4k+2$ case}

In this case $d_1, d_2: \eu_-\rightarrow \eu_-$
 and the inner product $C^{-1}$ when restricted  on $\eu_-$
 is non-degenerate. We define the adjoints $\delta_1, \delta_2:
 \eu_-\rightarrow \eu_-$
 of the $d_1, d_2$ operators, respectively. Again there are two Laplace operators
 $\Delta_1=\delta_2 d_1+d_1\delta_2$ and $\Delta_2=\delta_1 d_2+d_2\delta_1$
The expressions of these operators are the same as those in the previous section.
In particular we have that

\begin{corollary}

The $\Delta_1$ Laplace operator is
\bea
\Delta_1\phi&=&-(-1)^{s_C+s_\Gamma} C(\zeta_1, \zeta_2) \biggl(
 (g^{\mu\nu}+i\Omega^{\mu\nu})
\nabla_\mu \nabla_\nu \phi
\cr
&+& {1\over{\rm dim}_{\bC}(\Delta_+)}
\sum_{p=0}^{{m\over2}} {1\over (2p)!}
\cr
&&
\biggl[{i^{p+1}
(2p+2)!\over 2^{p+1} (2p+1) (p+1)!} \Omega^{\mu\nu} \Omega^{\rho_1\rho_2}\dots
\Omega^{\rho_{2p-1}\rho_{2p}} C\Gamma_{\rho_1\dots \rho_{2p}}
\rc R_{\mu\nu}\phi
\cr
&+& {i^{p-1} (2p)!\over 2^{p-2} (p-1)!} \Omega^{\rho_2\rho_3}\dots
\Omega^{\rho_{2p-2}\rho_{2p-1}} g^{\rho_1 \mu} g^{\rho_{2p} \nu} C\Gamma_{\rho_1\dots \rho_{2p}}
\rc R^{1,1}_{\mu\nu}\phi
\biggl]\biggr)~.
\eea
Similarly, the $\Delta_2$ Laplace operator is
\bea
\Delta_2\phi&=&-(-1)^{s_C+s_\Gamma} C(\zeta_1, \zeta_2) \biggl( (-1)^{s_c}
 (g^{\mu\nu}-i\Omega^{\mu\nu})
\nabla_\mu \nabla_\nu \phi
\cr
&+& {(-1)^{s_c}\over{\rm dim}_{\bC}(\Delta_+)}
\sum_{p=0}^{{m\over2}}
\biggl[{(-i)^{p+1}
(2p+2)!\over 2^{p+1} (2p+1) (p+1)!} \Omega^{\mu\nu} \Omega^{\rho_1\rho_2}\dots
\Omega^{\rho_{2p-1}\rho_{2p}} C\Gamma_{\rho_1\dots \rho_{2p}}
\rc R_{\mu\nu}\phi
\cr
&+& {(-i)^{p-1} (2p)!\over 2^{p-2} (p-1)!} \Omega^{\rho_2\rho_3}\dots
\Omega^{\rho_{2p-2}\rho_{2p-1}} g^{\rho_1 \mu} g^{\rho_{2p} \nu} C\Gamma_{\rho_1\dots \rho_{2p}}
\rc R^{1,1}_{\mu\nu}\phi
\biggl]\biggr)~.
\eea
\end{corollary}

One could also define two more Laplace operators $\hat\Delta_1=d_1\delta_1+\delta_1 d_1$
and $\hat\Delta_2=d_2\delta_2+\delta_2 d_2$. However they vanish. This can be
seen by a direct computation
\bea
\hat\Delta_1\phi&=&-(-1)^{s_C+s_\Gamma} \biggl( C( \zeta_1, \Gamma^\nu \Gamma^\mu  \zeta_1)
\nabla_\mu \nabla_\nu \phi
\cr
&+&
{1\over {\rm dim}_{\bC}(\Delta_+)} \sum_{p=0}^{{m\over2}} {1\over (2p)!}
\biggl[C(\zeta_1, \Gamma^{\mu\rho_1\dots \rho_{2p}\nu}  \zeta_1)
\cr
&+&(2p) (2p-1)g^{\mu \rho_1} C(\zeta_1, \Gamma^{\rho_2\dots \rho_{2p-1}}\zeta_1)
 g^{\rho_{2p}\nu}
\biggr]
 C\Gamma_{\rho_1\dots \rho_{2p}}
\rc R_{\mu\nu}\phi=0
~,
\eea
because the (0,2) part of the curvature $R$ vanishes.  Similarly $\hat\Delta_2=0$.

\newsection{$SU(m)$ holonomy, twisted complexes and algebraic spin cohomology}

\subsection{Twisted $\eu_\pm\otimes E$.}

Let $M$ a Riemannian manifold equipped with a connection $\nabla^M$
such that ${\rm hol}(\nabla^M)\subseteq SU(m)$. In addition let $E$ be a vector bundle
over $M$ equipped with a connection $\nabla^E$ and associated curvature $F$.

As in the previous section one can construct two first order
differential spin operators $d_1$ and $d_2$ associated with the
two parallel spinors of $\nabla^M$. Then one can use the vector bundle $E$
to twist the complexes $\eu_+$ and $\eu_-$ as $\eu_+\otimes E$ and $\eu_-\otimes E$,
respectively. Furthermore, one can
use the connection $\nabla^E$ to extend $d_1$ and $d_2$ to the twisted complexes
as it has been described in (\ref{twtwe}). We denote the
extended operators  with the same symbols.
 In particular we have

\begin{corollary}
The operators $d_1: \eu_+\otimes E\rightarrow \eu_+\otimes E$ for
$m=4k+1, 4k+3$ and $d_1: \eu_-\otimes E\rightarrow \eu_-\otimes E$ for
$m=4k, 4k+2$ are nilpotent, $d_1^2=0$, if the (2,0) part of the
curvature, $R$ and $F$,
of both the connections $\nabla$ and $\nabla^E$ vanishes. Similarly,
The operator $d_2: \eu_-\otimes E\rightarrow \eu_-\otimes E$, $m=4k, 4k+1, 4k+2, 4k+3$,
 is nilpotent,
 if the (0,2) part of the
curvature, $R$ and $F$,
of both the connections $\nabla$ and $\nabla^E$ vanishes. For $m=4k, 4k+2$,
$d_1 d_2+d_2 d_1=0$, if $F=R=0$.
\end{corollary}

The Laplace operators can be easily computed. In particular we find that
\bea
\Delta_1\phi&=&-(-1)^{s_C+s_\Gamma} C(\zeta_1, \zeta_2) \biggl(
 (g^{\mu\nu}+i\Omega^{\mu\nu})
\nabla_\mu \nabla_\nu \phi
\cr
&+& {1\over{\rm dim}_{\bC}(\Delta_+)}
\sum_{p=0}^{[{m\over2}]} {1\over (2p)!}
\cr
&&
\biggl[{i^{p+1}
(2p+2)!\over 2^{p+1} (2p+1) (p+1)!} \Omega^{\mu\nu} \Omega^{\rho_1\rho_2}\dots
\Omega^{\rho_{2p-1}\rho_{2p}} C\Gamma_{\rho_1\dots \rho_{2p}}
\rc (R_{\mu\nu}+F_{\mu\nu})\phi
\cr
&+& {i^{p-1} (2p)!\over 2^{p-2} (p-1)!} \Omega^{\rho_2\rho_3}\dots
\Omega^{\rho_{2p-2}\rho_{2p-1}} g^{\rho_1 \mu} g^{\rho_{2p} \nu} C\Gamma_{\rho_1\dots \rho_{2p}}
\rc (R^{1,1}_{\mu\nu}+F^{1,1}_{\mu\nu})\phi
\biggl]\biggr)~.
\la{laptone}
\eea
Similarly, the $\Delta_2$ Laplace operator is
\bea
\Delta_2\phi&=&-(-1)^{s_C+s_\Gamma} C(\zeta_1, \zeta_2) \biggl( (-1)^{s_c}
 (g^{\mu\nu}-i\Omega^{\mu\nu})
\nabla_\mu \nabla_\nu \phi
\cr
&+& {(-1)^{s_c}\over{\rm dim}_{\bC}(\Delta_+)}
\sum_{p=0}^{[{m\over2}]}
\biggl[{(-i)^{p+1}
(2p+2)!\over 2^{p+1} (2p+1) (p+1)!} \Omega^{\mu\nu} \Omega^{\rho_1\rho_2}\dots
\Omega^{\rho_{2p-1}\rho_{2p}} C\Gamma_{\rho_1\dots \rho_{2p}}
\rc (R_{\mu\nu}+F_{\mu\nu})\phi
\cr
&+& {(-i)^{p-1} (2p)!\over 2^{p-2} (p-1)!} \Omega^{\rho_2\rho_3}\dots
\Omega^{\rho_{2p-2}\rho_{2p-1}} g^{\rho_1 \mu} g^{\rho_{2p} \nu} C\Gamma_{\rho_1\dots \rho_{2p}}
\rc (R^{1,1}_{\mu\nu}+F^{1,1}_{\mu\nu})\phi
\biggl]\biggr)~.
\la{lapttwo}
\eea

\subsection{Twisted $\eu_\pm(E)$ complexes}

Suppose that the parallel spinors with respect to the ${\cal D}$ connection
on $\Delta_\pm \otimes E$ are in the direction of either
 $1\otimes 1$ or $e^1\wedge\dots\wedge e^m\otimes 1$
and $E$ is a vector bundle with a fibre inner product $h$, $\nabla^E h=0$.
One can construct an invariant inner product on $\Delta_\pm\otimes E$ as
$C^{-1}\otimes h$ and extended to the twisted complexes $\eu_\pm(E)$.
We can again define operators $d_1$ and $d_2$. In particular we have

\begin{proposition}
The operators  $d_1: \eu_+(E)\rightarrow \eu_+(E)$ for
$m=4k+1, 4k+3$ and $d_1: \eu_-(E)\rightarrow \eu_-(E)$ for
$m=4k, 4k+2$ are nilpotent, $d_1^2=0$,  if the (2,0) part of the
curvature of ${\cal D}$ vanishes. Similarly,
The operator $d_2: \eu_-(E)\rightarrow \eu_-(E)$,  $m=4k, 4k+1, 4k+2, 4k+3$,
is nilpotent, $d_2^2=0$, if the (0,2) part of the
curvature of ${\cal D}$ vanishes. For $m=4k, 4k+2$,
$d_1 d_2+d_2 d_1=0$, if the curvature of ${\cal D}$ vanishes.
\end{proposition}
\bproof
The proof is similar to that we have  already
investigated in the previous sections for the untwisted  $d_1$ and $d_2$ operators.
However, there is one difference.
If $E$ is not a line bundle, then in the expression
for $d_1^2$ and $d_2^2$ both symmetric and skew-symmetric
$C\Gamma^{(p)}$
contribute. This is unlike the untwisted case where only
the skew-symmetric $C\Gamma^{(p)}$ contribute.
However, there is no additional restriction on the
curvature of ${\cal D}$.
\eproof

The Laplace operators $\Delta_1$ and $\Delta_2$ can be easily computed in this case.
The expressions are as in (\ref{laptone}) and (\ref{lapttwo}) with the curvature
$R$ and $F$ replaced by the curvature of ${\cal D}$.

\subsection{Algebraic Cohomologies}

We have seen that the operator $D$ on the complex $\Lambda^*\otimes\eu$
is defined provided that
$C\Gamma$ is skew-symmetric which is the case for $C=A$ if $m=4k, 4k+3$
 and for $C=B$ if $m=4k+2$.
Moreover $D$ restricts on $\eu_\pm$, $D: \eu_\pm\rightarrow \eu_\pm$
if $m=4k+3$. Therefore we conclude that the $D$
 operator can be defined on the complexes
$\Lambda^*\otimes \eu_\pm$  and $\Lambda^*\otimes \eu_\pm\otimes E$ only for $m=4k+3$.
The operator $D$ can also defined  for twisted complexes $\Lambda^*\otimes \eu_\pm(E)$
but we shall not investigate this further here.

\begin{corollary}
The algebraic operator $D$ anticommutes with both $d_1$
and $d_2$ differential operators
\be
d_1D+D d_1=d_2 D+ D d_2=0~.
\ee
\end{corollary}

Therefore one can define the operators $d_1+D$ and $d_2+D$ which are nilpotent
provided $d_1^2=d_2^2=0$. The cohomology of $d_1+D$ and $d_2+D$ can be computed
using spectral sequences. We shall not do a general computation. Instead,
we shall give the cohomology of the operator $d_2+D$ in the special case where
$M$ is a six-dimensional Calabi-Yau manifold.

\newsection{Complex manifolds with holonomy $SU(m)$ and spin cohomology}

It is clear from the results of the previous section that complex  spin
cohomology  is related
to the Dolbeault cohomology. Here we shall establish the precise relation
and we shall give the classes of the spin cohomology in terms of those
of the Dolbeault cohomology, see e.g. \cite{chern, harris}.

\subsection{ Spin and Dolbeault cohomologies}

Let $M$ be a complex manifold equipped with a connection $\nabla$, ${\rm hol}(\nabla)\subseteq SU(m)$.
On $M$, it is known that
\be
\Delta=\oplus_q\Lambda^{0,q}=\Lambda^{0,*}~.
\ee
This can be easily seen from $\Delta=\Lambda^{0,q}(1)$, where $\Lambda^{0,q}$ acts on $1$ with Clifford multiplication.
In particular, we have $\Delta_+=\Lambda^{0,{\rm even}}$ and $\Delta_-=\Lambda^{0,{\rm odd}}$.
Thus
\bea
\eu_+&=&\Lambda^*(\Lambda^{0,{\rm even}})
\cr
\eu_-&=&\Lambda^*(\Lambda^{0,{\rm odd}})
\eea
Write $\Delta_-=\Lambda^{0,1}\oplus Z$, where $Z=\oplus_{p\geq 1}\Lambda^{0,2p+1}$.
The complex $\eu_-$ can now be decomposed as
\be
\eu_-^\ell=\oplus_{p+q=\ell} \Lambda^{0,p}\otimes \Lambda^q(Z)~.
\ee

\begin{proposition}
\be
d_2=\bar \partial: \Lambda^{0,p}\otimes \Lambda^q(Z)\rightarrow
\Lambda^{0,p+1}\otimes \Lambda^q(Z)~.
\ee
\end{proposition}
\bproof
To show this, we first evaluate the action of $d_2$ on $\Lambda^{0,1}\subset \eu_-$. Indeed let
$\eta_i e^i\in \Lambda^{0,1}$, then  we have
\be
d_2 (\eta_i e^i)= \lambda (-1)^{{1\over2} m(m-1)} (\nabla_i+i\nabla_{i+m}) \eta_j~ e^i\bar\wedge e^j=
\lambda (-1)^{{1\over2} m(m-1)} \bar\partial_i \eta_j~ e^i\wedge e^j~,
\ee
where $\lambda=1$ for the $A$ inner product and $\lambda=i^m$ for $B$
inner product. After suppressing the numerical coefficient which is inconsequential for the computation of cohomology,
we have $d_2=\bar\partial$ on $\Lambda^{0,1}$. Using the definition of the $\bar\wedge$ product, it is straightforward
to extend the proof to the rest of the complex $\eu_-$.
\eproof

\begin{corollary}
Let $M$ be a complex manifold as described in the beginning of the section. Then the spin cohomology
\be
H^\ell_{d_2}(\eu_-)=\oplus_{p+q=\ell}H^{0,p}_{\bar\partial}(\Lambda^q(Z))~.
\ee
Therefore the spin cohomology of the $d_2$ operator can be computed in terms
of  Dolbeault cohomology of a twisted complex by the bundle $\Lambda^*(Z)$, where $Z=\oplus_{p\geq 1}\Lambda^{0,2p+1}$.
\end{corollary}

A direct consequence of this is that
\be
H^\ell_{d_2}(\eu_-\otimes E)=\oplus_{p+q=\ell}H^{0,p}_{\bar\partial}(\Lambda^q(Z)\otimes E)~.
\ee
Using the corollary, we can also compute the index of the spin complex $(\eu_-, d_2)$
 in terms of the index of the twisted
$\bar \partial$ complex. In particular, we have
\be
{\rm Index}_{d_2} (\eu_-)=\sum_{q\geq 0} (-1)^q {\rm Index}_{\bar\partial} (\Lambda^q(Z))
\ee
or more generally
\be
{\rm Index}_{d_2} (\eu_-\otimes E)=\sum_{q\geq 0} (-1)^q {\rm Index}_{\bar\partial} (\Lambda^q(Z)
\otimes E)~.
\ee

It remains to investigate the cohomology of $d_1$. We shall consider the cases
$m=4k+1, 4k+3$ and $m=4k, 4k+2$ separately.
In the former case $d_1: \eu_+ \rightarrow \eu_+$. Writing $\Delta_+=\Lambda^{0,m-1}\oplus W$,
where $W=\oplus_{p<{m-1\over2}}\Lambda^{0,2p}$, we
have
\be
\eu^\ell_+=\oplus_{p+q=\ell} \Lambda^p(\Lambda^{0,m-1})\otimes \Lambda^q(W)~,
\ee
and
\be
d_1: \Lambda^p(\Lambda^{0,m-1})\otimes \Lambda^q(W)\rightarrow
\Lambda^{p+1}(\Lambda^{0,m-1})\otimes \Lambda^q(W)~.
\ee
Since there is a $SU(m)$ structure, we can identify $\Lambda^{0,m-1}=\Lambda^{1,0}$,
$\Lambda^p(\Lambda^{0,m-1})=\Lambda^{p,0}$ and
\be
d_1=\partial: \Lambda^{p,0}\otimes\Lambda^q(W)\rightarrow
\Lambda^{p+1,0}\otimes \Lambda^q(W)~.
\ee
Therefore, we conclude that
\be
H^\ell_{d_1}(\eu_+)=\oplus_{p+q=\ell}H^{p,0}_{\partial}(\Lambda^q(W))~.
\ee
Thus $H^*_{d_1}(\eu_+)=H^*_{d_2}(\eu_-)$. The same applies for
$m=4k, 4k+2$.

\subsection{Algebraic operations on twisted complexes}


On complex manifolds with ${\rm hol}(\nabla)\subseteq SU(m)$ apart from the
differential spin operators $d_1, d_2$,
there is also an algebraic spin operator $D$ provided
$m=4k+3$. We shall  focus on the twisted spin cohomology
associated with the operators $d_2$ and $D$. In the context
of complex geometry, there are several versions that one can
consider. In particular one can defined a twisted spin cohomology
on $\Lambda^*\otimes \eu_-$ as we have already mentioned in section 5.3.
However it is also possible to twist $\eu_-$ with either
$\Lambda^{*,0}$ or $\Lambda^{0,*}$. In each of these cases the
twisted spin cohomology of the operator $d_2+D$, or equivalently $\bar\partial+D$, can be
computed using a spectral sequence.

The twisted complex $C=\Lambda^*\otimes \eu_-$ is a double complex, $C^{p, \ell}$, with grading induced
from the space of forms $\Lambda^*$
and that of $\eu_-$. However in this grading, $d_2=\bar\partial$ and $D$ do not act
with horizontal and vertical operations. In particular, $d_2: C^{p, \ell}\rightarrow C^{p, \ell+1}$ and
$D: C^{p, \ell}\rightarrow C^{q-1, \ell+2}$. It is therefore convenient
to introduce a new grading  as
\be
C^{[-p, \ell+2p]}= C^{p, \ell}=\Lambda^p\otimes \eu^\ell_-~.
\ee
Note know that $d_1:  C^{[-p, \ell+2p]}\rightarrow C^{[-p, \ell+2p+1]}$
and $D: C^{[-p, \ell+2p]}\rightarrow C^{[-p+1, \ell+2p]}$ as expected.
The twisted complexes $\Lambda^{*,0}\otimes \eu_-$ and $\Lambda^{0,*}\otimes \eu_-$
 can be treated in a similar way.
The machinery of spectral sequences can now be used to do the computation, see e.g \cite{bott} and references within.
Instead of developing the general theory of computing the cohomology of $\bar\partial+D$
for the various complexes above,
 we shall give the cohomology of $(\Lambda^{*,0}\otimes \eu_-, \bar\partial+D)$
 for six-dimensional Calabi-Yau manifolds in an example below.

\subsection{$Spin_c$ structures and spin cohomology}

Let $M$ be a complex  manifold equipped with a $Spin_c$ structure and compatible connection ${\cal \nabla}$,
${\rm hol}({\cal\nabla})\subseteq SU(m)$. Suppose that $L$ is a (locally defined) complex line bundle
associated with the $Spin_c$ structure.
On $M$, it is known that
\be
\Delta^*\otimes L=\oplus_q\Lambda^{0,q}=\Lambda^{0,*}~.
\ee
This is similar to the standard complex case, we have investigated.
In particular, we have $\Delta_+\otimes L=\Lambda^{0,{\rm even}}$ and $\Delta_-=\Lambda^{0,{\rm odd}}$.
Thus
\bea
\eu_+(L)&=&\Lambda^*(\Lambda^{0,{\rm even}})
\cr
\eu_-(L)&=&\Lambda^*(\Lambda^{0,{\rm odd}})
\eea
Write $\Delta_-\otimes L=\Lambda^{0,1}\otimes L\oplus Z\otimes L$, where $Z=\oplus_{p\geq 1}\Lambda^{0,2p+1}$.
The complex $\eu_-(L)$ can now be decomposed as
\be
\eu_-^\ell(L)=\oplus_{p+q=\ell} \Lambda^{p}(\Lambda^{0,1})\otimes \Lambda^q(Z)~.
\ee

\begin{proposition}
\be
d_2=\bar \partial: \Lambda^{p}(\Lambda^{0,1})\otimes \Lambda^q(Z\otimes L)\rightarrow
\Lambda^{p+1}(\Lambda^{0,1})\otimes \Lambda^q(Z)~.
\ee
\end{proposition}
\bproof
To show this, we first evaluate the action of $d_2$ on $\Lambda^{0,1}\subset \eu_-(L)$. Indeed let
$\eta_i e^i\in \Lambda^{0,1}$, then  we have
\be
d_2 (\eta_i e^i)= \lambda (-1)^{{1\over2} m(m-1)} ({\cal\nabla}_i+i{\cal \nabla}_{i+m}) \eta_j~ e^i\bar\wedge e^j=
\lambda (-1)^{{1\over2} m(m-1)} \bar\partial_i \eta_j~ e^i\wedge e^j~,
\ee
where $\lambda=1$ for the $A$ inner product and $\lambda=i^m$ for $B$
inner product. After suppressing the numerical coefficient which is inconsequential for the computation of cohomology,
we have $d_2=\bar\partial$ on $\Lambda^{0,1}$. Using the definition of the $\bar\wedge$ product, it is straightforward
to extend the proof to the rest of the complex $\eu_-(L)$.
\eproof

\begin{corollary}
Let $M$ be a complex $Spin_c$ manifold as described in the beginning of the section. Then the spin cohomology
\be
H^\ell_{d_2}(\eu_-(L))=\oplus_{p+q=\ell}H^{0,p}_{\bar\partial}(\Lambda^q(Z))~.
\ee
Therefore the spin cohomology of the $d_2$ operator can be computed in terms
of  Dolbeault cohomology of a twisted complex by the bundle $L^p\otimes\Lambda^q(Z)$,
 where $Z=\oplus_{p\geq 1}\Lambda^{0,2p+1}$.
\end{corollary}

\newsection{Spin cohomology and six-dimensional Calabi-Yau manifolds}

\subsection{Differential spin cohomology}

Applying the general theory of the previous section to this case, we have
$\Delta_-=\Lambda^{0,{\rm odd}}=\Lambda^{0,1}\oplus \Lambda^{0,3}$. In addition
for six-dimensional Calabi-Yau manifolds
$\Lambda^{0,3}$ is trivial line bundle.
Using these, we find that
\be
\eu^\ell_-=\Lambda^{0,\ell}\oplus \Lambda^{0,\ell-1}~
\la{eul}
\ee
and
\be
d_2=\bar \partial: \Lambda^{0,\ell}\oplus \Lambda^{0,\ell-1}\rightarrow
\Lambda^{0,\ell+1}\oplus \Lambda^{0,\ell}~.
\ee
Therefore
\be
H_{d_2}^\ell(\eu_-)=H^{0,\ell}_{\bar\partial}\oplus H^{0,\ell-1}_{\bar\partial}~.
\ee
In particular for an irreducible six-dimensional Calabi-Yau manifold, we have that
\be
H^0_{d_2}=\bC~,~~H^1_{d_2}=\bC~,~~H^2_{d_2}=0~,~~H^3_{d_2}=\bC~,~~H^4_{d_2}=\bC~.
\ee

It is also straightforward to compute the cohomology of the twisted complex
$(\Lambda^{*,0}\otimes \eu_-, d_2)$. In particular, we find that
\be
H_{d_2}^{p,\ell}(\Lambda^{*,0}\otimes\eu_-)=H^{p,\ell}_{\bar\partial}\oplus H^{p,\ell-1}_{\bar\partial}~.
\ee

\subsection{Twisted complexes and algebraic spin cohomology}

To compute the cohomology of the complex $(\Lambda^{*,0}\otimes \eu_-, d_2+D)$, we first investigate the
complex $(\Lambda^{*,0}\otimes \eu_-, D)$.
Using  (\ref{eul}), we find that the operator $D: \Lambda^{p,0}\otimes \eu^\ell_-\rightarrow \Lambda^{p-1,0}\otimes \eu^{\ell+2}_- $
acts  as
\be
D:  \Lambda^{p,0}\otimes [\Lambda^{0,\ell}\oplus \Lambda^{0,\ell-1}]\rightarrow
\Lambda^{p-1}\otimes [\Lambda^{0,\ell+2}\oplus \Lambda^{0,\ell+1}]~
\ee
or equivalently
$D: \Lambda^{p,\ell}\oplus \Lambda^{p,\ell-1}\rightarrow \Lambda^{p-1, \ell+2}\oplus \Lambda^{p-1,\ell+1}$.
Since it acts on the two parts in the sum separately, it is enough to consider
only its action in the first part.
After some computation, one finds that
\be
D\psi=- (-1)^q{1\over (p-1)! q! 2}\psi_{\gamma
\alpha_1\dots \alpha_{p-1}, \bar\beta_3\dots\bar\beta_{q+2}} \epsilon^\gamma{}_{\bar\beta_1\bar\beta_2} e^{\alpha_1}\wedge\dots
\wedge e^{\alpha_{p-1}}\wedge e^{\bar\beta_1}\wedge e^{\bar\beta_2}\wedge\dots
\wedge e^{\bar\beta_{q+2}}~,
\ee
where $\psi\in \Lambda^{p,q}$.

Next we shall compute the cohomology of the double complex $(\Lambda^{*,*}, \bar\partial+D)$
 using a spectral sequence, see e.g. \cite{bott} and references within.
The most convenient filtration is that for which
$E_1=H^{p,q}_{\bar \partial}$. Then from the general theory of spectral sequences for
double complexes $E_2=H_{D} H_{\bar\partial}$ and $E_2$ is graded as
the double complex.
It is known that for six-dimensional irreducible Calabi-Yau manifolds the non-vanishing Dolbeault groups
are  $H^{0,0}_{\bar \partial}$, $H^{1,1}_{\bar \partial}$, $H^{2,1}_{\bar \partial}$, $H^{1,2}_{\bar \partial}$
$H^{2,2}_{\bar \partial}$ and $H^{3,3}_{\bar\partial}$. Moreover $H^{3,0}_{\bar \partial}=H^{0,3}_{\bar \partial}=\bC$ which are
generated by the parallel (3,0)- and (0,3)-forms, respectively. To compute $E_2$ observe that
\be
0\rightarrow H^{3,0}_{\bar \partial}{\buildrel D\over \rightarrow} H^{2,2}_{\bar \partial}\rightarrow 0
\ee
and
\be
0\rightarrow H^{1,1}_{\bar \partial}{\buildrel D\over \rightarrow} H^{0,3}_{\bar \partial}\rightarrow 0~,
\ee
and the rest of the cohomology groups of $E_1=H^{p,q}_{\bar \partial}$ live in $E_2$.
It is easy to see that ${\rm Ker}D|_{H^{3,0}_{\bar \partial}}=\{0\}$. Therefore
$E_2^{3,0}$ vanishes. In addition $D (H^{3,0}_{\bar \partial})=\bC<\Omega\wedge\Omega>$, where $\Omega$ is the
K\"ahler form. Since $D H^{2,2}_{\bar \partial}=0$,
we conclude that $E_2^{2,2}=H^{2,2}_{\bar \partial}/\bC<\Omega\wedge\Omega>$, ie $E_2^{2,2}=PH^{2,2}_{\bar \partial}$
is generated by the
primitive $(2,2)$ harmonic forms.

Next observe that
\be
{\rm Ker}D|_{H^{1,1}_{\bar \partial}}=\{\alpha \in H^{1,1}_{\bar \partial}~~{\rm such~that}~~\Omega\cdot\alpha=0 \}~.
\ee
Therefore $E_2^{1,1}=P{H^{1,1}_{\bar \partial}}$ is generated by the primitive $(1,1)$ harmonic forms. In addition we have that
$D H^{1,1}_{\bar \partial}= H^{0,3}_{\bar \partial}$, therefore $E_2^{0,3}=0$. Thus the only non-vanishing groups
are $E_2^{0,0}=\bC, E_2^{2,1}=H^{2,1}_{\bar\partial}, E_2^{1,2}=H^{1,2}_{\bar\partial}, E^{1,1}_2 =P{H^{1,1}_{\bar \partial}},
E^{2,2}_2 =P{H^{2,2}_{\bar \partial}}$ and $E^{3,3}_2=\bC$.

It remains to show that $E_2=E_\infty$. This
is easily verified by computing the action of the differential $d_2$  of $E_2$.
For this, we need  to convert to the grading  of the   double complex  $\Lambda^{*,*}$
for which $\bar\partial: \Lambda^{[m,n]}\rightarrow \Lambda^{[m,n+1]}$ acts vertically
and $D: \Lambda^{[m,n]}\rightarrow \Lambda^{[m+1,n]}$ acts horizontally. As we have explained
$(m,n)=(-p, 2p+q)$, ie $E_2^{[-p, 2p+q]}= H_D^p H^q_{\bar\partial}$. The differential
$d_2: E_2^{[m,n]}\rightarrow E_2^{[m-2, n-1]}$.
It is easy then to see  that the $d_2$ differential is the zero map and $E_2=E_\infty$. Therefore the
cohomology of the operator $d+D$ is given by $E_2$. Thus we have shown the proposition.

\begin{proposition}
Let $M$ be an irreducible six-dimensional Calabi-Yau manifold. The non-vanishing cohomology groups of the complex $(\Lambda^{p,q}, \bar\partial+D)$ are
 $H^0_{\bar\partial+D}=H^{0,0}_{\bar\partial}=\bC, H^2_{\bar\partial+D}=P{H^{1,1}_{\bar \partial}},
 H^3_{\bar\partial+D}=H^{2,1}_{\bar\partial}\oplus H^{1,2}_{\bar\partial},  H^4_{\bar\partial+D}=P{H^{2,2}_{\bar \partial}}$
 and $H^6_{\bar\partial+D}=H^{3,3}_{\bar\partial}=\bC$.
\end{proposition}

Finally, we have $H^\ell_{\bar\partial+D}(\Lambda^{*,0}\otimes \eu_-)= H^\ell_{\bar\partial+D}\oplus H^{\ell-1}_{\bar\partial+D}$.

\newsection{Manifolds with connections with holonomy $Sp(k)$ and spin cohomology}

Let $M$ be a Riemannian manifold which admits a connection $\nabla$ with holonomy $Sp(k)$.
In this case, there are $k+1$ parallel spinors and so $k+1$ spin differential operators
that one can define. The spin differential operators associated with the parallel
spinors $1$ and $e_1\wedge\dots\wedge e_{2k}$ are the same as the $d_1$ and $d_2$
spin differential operators that we have investigated for complex manifolds. There are
another $k-1$ spin differential operators $d_{\omega^k}$ associated with the $\omega^k$, $1\leq k\leq k-1$, parallel spinors
of section two. Since $\omega\in \Lambda^2$ and $m=2k$ even, $d_{\omega^k}: \eu_-\rightarrow \eu_-$.
We shall not present a complete analysis of all  cases. Instead we shall
focus on the spin differential operator $d_0=d_{\omega^{k-1}}$ associated with the parallel spinor
$\omega^{k-1}$.

As we have explained for complex manifolds, $\eu_-$ can  be decomposed as
\be
\eu_-^\ell=\oplus_{p+q=\ell} \Lambda^{0,p}\otimes \Lambda^q(Z)~.
\ee

\begin{proposition}
Let $M$ be a hyperK\"ahler manifold, then
\be
d_0=K\lc \partial: \Lambda^{0,p}\otimes \Lambda^q(Z)\rightarrow
\Lambda^{0,p+1}\otimes \Lambda^q(Z)~,
\ee
where $K$ is the second complex structure on $M$.
\end{proposition}
\bproof
To show this, we first evaluate the action of $d_0$ on $\Lambda^{0,1}\subset \eu_-$. Indeed let
$\eta_i e^i\in \Lambda^{0,1}$, then  we have
\be
d_0 (\eta_i e^i)= \lambda <\omega, \Gamma^i e_k> (\nabla_i-i\nabla_{i+m}) \eta_j~ e^k\bar\wedge e^j=
\lambda  K^i{}_k\partial_i \eta_j~ e^k\wedge e^j~,
\ee
where $\lambda$ is inconsequential  numerical coefficient that depends on the inner product and the
normalization of the parallel spinor and $<v, Kw>=\omega(v,w)$. Thus we have
that $d_2=K\lc \partial$ on $\Lambda^{0,1}$. Using the definition of the $\bar\wedge$ product, it is straightforward
to extend the proof to the rest of the complex $\eu_-$.
\eproof

\newsection{Real Spin Cohomologies}

So far we have investigated the complex spinor cohomologies.
 Now we shall turn to investigate the
real ones. The real spinor representations for $m=4k, 4k+1, 4k+3$ can be constructed
by imposing a reality condition on the complex representations.
These reality conditions
are
\bea
\eta&=&\pm A(\bar\eta)~,~~~~~~\eta\in \Delta^\pm~,~~~~~~~~~~~~~~m=4k
\cr
\eta&=&A(\bar\eta)~,~~~~~~~~\eta\in \Delta^+\oplus\Delta^-~,~~~~=4k+1
\cr
\eta&=&B(\bar\eta)~,~~~~~~~~\eta\in \Delta^+\oplus\Delta^-~,~~~~=4k+3
\eea
and $\Delta_{\bR}^\pm=\{\eta\in \Delta^\pm,~{\rm s.t.}~\eta=\pm A(\bar\eta)\}$,
 $\Delta_{\bR}=\{\eta\in \Delta^+\oplus\Delta^-,~{\rm s.t.}~\eta=A(\bar\eta)\}$
 and  $\Delta_{\bR}=\{\eta\in \Delta^+\oplus\Delta^-,~{\rm s.t.}~\eta=B(\bar\eta)\}$,
 respectively.

\subsection{$SU(m)$ invariant spinors}

We begin with a summary of the properties of real spin representation in various cases.

\subsubsection{m=4k}

The real parallel spinors in the $SU(4k)$ case are
\be
\tau_1={1\over \sqrt{2}} (1+e_1\wedge\dots\wedge e_{m})~,~~\tau_2={i\over \sqrt{2}} (1-e_1\wedge\dots\wedge e_m)~.
\ee
Both parallel spinors $\tau_1, \tau_2\in \Delta_{\bR}^+$. We can again define
spin cohomologies $d_1, d_2$  associated with $\tau_1, \tau_2$ for $C=A$ on
the real complex $\eu^{\bR}_{-}$.
Thus $s_c=0, s_\Gamma=1$. The decomposition of the real spinor representations
is
\be
\Delta_{\bR}^\pm\otimes \Delta_{\bR}^\pm=\sum^{2k-1}_{p=0} \Lambda_{\bR}^{2p}\oplus
\Lambda^{m\pm}_{\bR}~.
\ee

\subsection{m=4k+1}

The real parallel spinors are
\be
\tau_1={1\over \sqrt{2}} (1+e_1\wedge\dots\wedge e_{m})~,~~
\tau_2={i\over \sqrt{2}} (1-e_1\wedge\dots\wedge e_m)~.
\ee
The parallel spinors $\tau_1, \tau_2\in \Delta_{\bR}$. We can again define
spin cohomologies $d_1, d_2$  associated with $\tau_1, \tau_2$ for $C=A$ on
the real complex $\eu^{\bR})$.
Thus $s_c=0, s_\Gamma=0$. The decomposition of the real spinor representations
is
\be
\Delta_{\bR}\otimes \Delta_{\bR}=\sum^{2m}_{p=0} \Lambda_{\bR}^{2p}~.
\ee

\subsection{m=4k+3}

The real parallel spinors are
\bea
\tau_1&=&{1\over \sqrt{2}} (1+i e_1\wedge\dots\wedge e_{m})
\cr
\tau_2&=&{1\over \sqrt{2}} (i 1+e_1\wedge\dots\wedge e_m)~.
\eea
The parallel spinors $\tau_1, \tau_2\in \Delta_{\bR}$. We can again define
spin cohomologies $d_1, d_2$  associated with $\tau_1, \tau_2$ for $C=B$ on
the real complex $\eu^{\bR}$.
Thus $s_c=0, s_\Gamma=1$. The decomposition of the real spinor representations
is
\be
\Delta_{\bR}\otimes \Delta_{\bR}=\sum^{2m}_{p=0} \Lambda_{\bR}^{2p}~.
\ee

In all the above cases we find the following:
\begin{theorem}
The operators $d_1, d_2$ are nilpotent iff the curvature of the connection
$\nabla$ vanishes.
\end{theorem}
\bproof
This is a consequence of the results we have already demonstrated in sections four and six.
\eproof

We also have that
\begin{theorem}
The Laplacians $\Delta_1, \Delta_2$ of the operators $d_1, d_2$ are
\bea
\Delta_2\phi=\Delta_1\phi&=& g^{\mu\nu} \nabla_\mu\nabla_\nu \phi
\eea
\end{theorem}
\bproof
This is a consequence of the results we have already demonstrated in section four.
\eproof

\begin{corollary}
The  real spin cohomologies  $H^*(\eu^{\bR})$ defined above are generated by the parallel elements in
$\eu^{\bR}$ with respect to $\nabla$.
\end{corollary}
\bproof
This follows from a partial integration formula and the fact that the inner
product  $C=A,B$ restricted in $\Delta_{\bR}$  is definite.
\eproof

A class of manifolds which we can define a real spin cohomology are group
manifolds equipped with the left or right invariant connections. One can also
defined twisted real spin cohomology but we shall not pursue this further here.

\newsection{$Spin(7)$ spin cohomology}

The $Spin(7)$ invariant spinor is $\zeta={1\over\sqrt{2}} (e_1-e_2\wedge e_3\wedge e_4)$.
Therefore $\zeta\in \Delta^-$. So the spin cohomology operator is
$d:\eu_+\rightarrow \eu_+$.

\begin{theorem}
The spin operator is nilpotent, $d^2=0$, if the connection $\nabla$ is the
Levi-Civita connection of a metric on the manifold $M$ with
holonomy $Spin(7)$.
\end{theorem}

\bproof
The representations $\Delta^\pm$ are real.
The map $\tau:  \Delta^+\rightarrow \Lambda^1(\bR^8)$ given by
$\tau(\eta)=C\Gamma_\mu(\zeta, \eta) e^\mu=C_\zeta^\mu(\eta) e_\mu$
induces an isomorphism between the $\Delta^+$ and the vector representations, $C=A$.
 This can be easily seen by
observing that there is a (real) basis in $\Delta^+$ such
that $\tau$ is diagonal. This basis is
\bea
&&1+e_1\wedge\dots\wedge e_4~,~~~i(1-e_1\wedge\dots\wedge e_4)~,~~~i(e_1\wedge e_2+e_3\wedge e_4)~,
\cr
&&(e_1\wedge e_2-e_3\wedge e_4)~,~~~e_1\wedge e_3+e_2\wedge e_4~,~~~~
i(e_1\wedge e_3-e_2\wedge e_4)~,
\cr
&&i(e_2\wedge e_3+e_1\wedge e_4)~,~~~(e_2\wedge e_3-e_1\wedge e_4)~.
\eea
In addition we have that
\be
d^2\phi={1\over2}C_\zeta^\mu\bar\wedge C_\zeta^\nu R_{\mu\nu}\phi~.
\ee
The right-hand-side will vanish if the curvature $R$ of the connection $\nabla$
is that of a Levi-Civita connection for a metric with holonomy $Spin(7)$ by virtue
of the Bianchi identity.
\eproof

\begin{corollary}
Let $M$ be a manifold equipped with a metric with holonomy $Spin(7)$. Then
$H^*(\eu_+)=H^*_{dR}(M)$.
\end{corollary}
\bproof
The map $\tau$ is an isomorphism between the
spin cohomology complex $(\eu_-,d)$ and the  de Rham complex $(\Lambda^*(M), d)$.
Therefore it induces an isomorphism in cohomology.
\eproof

\newsection{$G_2$ Spin Cohomology}

The spinor $Spin(7)$ representation $\Delta$ decomposes under $G_2$ as
$\Delta=\bR\oplus \Lambda^1(\bR^7)$. The $G_2$ invariant spinor is
$\zeta={1\over \sqrt{2}} (e_1-e_{234})$.
As in the previous cases, one can define a linear operator $d$ on $\eu$
using the spinor $\zeta$.

\begin{theorem}
Let $M$ a manifold equipped with a metric $g$ with holonomy contained in $G_2$.
The operator $d$ associated to the Levi-Civita connection is nilpotent.
\end{theorem}
\bproof
Since $\Delta=\bR\oplus \Lambda^1(\bR^7)$, we have that
$\eu^\ell=\Lambda^\ell(M)\oplus \Lambda^{\ell-1}(M)$. Moreover,
$\tau: \Delta\rightarrow \Lambda^1(\bR^7)$ such that $\tau(\eta)=C\Gamma^\mu(\zeta, \eta) e_\mu$
is onto and has kernel $\bR<\zeta>$. Next
\be
d^2\phi={1\over2} C\Gamma^\mu\bar\wedge C\Gamma^\nu \bar\wedge R_{\mu\nu}\phi~
\ee
which vanishes because of the Bianchi identity where $R$ is the curvature
of the $G_2$ metric.
\eproof

\begin{corollary}
Let $M$ a manifold equipped with a metric with holonomy $G_2$. Then
$H^\ell(\eu)=H^{\ell}_{dR}(M)\oplus H^{\ell-1}_{dR}(M)$.
\end{corollary}
\bproof
The map $\tau$ induces an isomorphism between the two complexes. This
establishes the isomorphism between the cohomologies.
\eproof

\subsection*{Acknowledgements}

I would like to thank P.S. Howe for many helpful discussions and for suggesting the term Spin Cohomology.
I would like to also thank Stefan Ivanov who took part in the initial stage of this project.

\end{document}